\documentclass[amsmath,amssymb,wasysym]{amsart}

\usepackage{graphicx}
\usepackage{amsmath}
\usepackage{amssymb}
\usepackage{amsfonts}
\usepackage{tikz}
\usetikzlibrary{tqft}
\usetikzlibrary{shapes}
\usepackage{cite}
\usepackage{tikz}
\usetikzlibrary{tqft}
\usetikzlibrary{shapes}
\usepackage{graphics}
\usepackage{float}
\usepackage{mathrsfs}
\usepackage{blindtext}
\usetikzlibrary{matrix,arrows}
\usepackage{tikz-cd}
\usepackage{caption}
\usepackage{pgfplots}
\usepgfplotslibrary{patchplots}
\usepackage[textwidth=6.2in, textheight=7.8in, hcentering]{geometry}

 \newtheorem{theorem}{Theorem}[subsection]
 
 \newtheorem{lemma}[theorem]{Lemma}
 \newtheorem{proposition}[theorem]{Proposition}
 \newtheorem{remark}[theorem]{Remark}
 
 \theoremstyle{definition}
 \newtheorem{definition}[theorem]{Definition}
 \theoremstyle{remark}
 \numberwithin{equation}{subsection}

\begin{document}

\title[Multiplicativity of Reidemister-Franz Torsion for Even Manifolds]{Multiplicativity of Reidemister-Franz Torsion for Even Manifolds}

\author{ESMA D\.{I}R\.{I}CAN ERDAL}\thanks{The author was partially supported by TUBİTAK under the project number 124F247.}
\address{}



\date{\today}

\dedicatory{}

\commby{}


\begin{abstract}
We study Reidemeister-Franz torsion for non-acyclic cellular chain complexes arising from closed, oriented, highly connected even dimensional manifolds. The monoid of such manifolds under connected sum admits a unique factorisation into indecomposable elements. Using this factorisation, we prove that the  Reidemeister-Franz torsion of an even-dimensional manifold decomposes multiplicatively as the product of the torsions of its prime factors without any corrective term. 
\end{abstract}


\maketitle

\section{Introduction}
\label{Sec:1}
All manifolds considered in this paper are assumed to be non-empty, closed, connected and oriented. An $n$-manifold $M^{n}$ is called highly connected if 
$$\pi_i(M^{n})=0$$ for $i=0,\ldots,\lfloor n/2\rfloor -1.$ We denote by  $\mathcal{M}_{n}^\mathrm{Diff}$ the set of the diffeomorphism classes of $n$-dimensional differentiable manifolds, and by $\mathcal{M}_{n}^{\mathrm{Diff},\mathrm{hc}}$ the subset consisting of highly connected manifolds.

	For $n\in \mathbb{N}$, let $M^{n}, N^{n} \in \mathcal{M}_{n}^{\mathrm{Diff}}.$ Given an orientation-preserving smooth embedding $\varphi:\overline{\mathbb{D}^{n}}\rightarrow M^n$ and an orientation-reversing smooth embedding $\varrho:\overline{\mathbb{D}^{n}}\rightarrow N^n,$ the connected sum of $M^{n}$ and $N^{n}$ is defined by
	 $$ M^{n} \# N^{n} = ((M^{n} \setminus \varphi({\mathbb{D}^{n}})\sqcup(N^{n} \setminus \varrho({\mathbb{D}^{n}}))/\sim,$$ 
	 where the equivalence relation $\sim$ identifies $\varphi(x)\sim \varrho(x) \mathrm{ \ for \  all \ } x \in \mathbb{S}^{n-1}=\partial\overline{\mathbb{D}^{n}}.$ The diffeomorphism type of the connected sum of two differentiable manifolds is independent of the choice of embedding, \cite[Theorem 2.7.4]{Wall2}. Hence, by \cite{Bokor2021}, 
$\mathcal{M}_{n}^{\mathrm{Diff}}$ and its subset $\mathcal{M}_{n}^{\mathrm{Diff},\mathrm{hc}}$ are abelian monoids (written multiplicatively) under connected sum operation.

\begin{definition} Let $\mathcal{M}$ be a monoid.
\begin{itemize}
\item[(i)] If $\mathcal{M}$ is abelian (written multiplicatively), then $m \in \mathcal{M}$ is called prime if it is not a unit and if it divides a product only if it divides one of the factors.
\item[(ii)] $\mathcal{M}^*$ denotes the units of 
$\mathcal{M}$ and we write $\bar{\mathcal{M}}:= \mathcal{M} /\mathcal{M}^*.$
\item[(iii)] If $\mathcal{M}$ is abelian, then $\mathcal{P}(\mathcal{M})$ denotes the set of prime elements in $\bar{\mathcal{M}}.$ Moreover, $\mathcal{M}$ is called a unique factorisation monoid if the canonical monoid morphism $\mathbb{N}^{\mathcal{P}(\mathcal{M})}\rightarrow \mathcal{M}$ is an isomorphism.
\end{itemize}
 \end{definition}

\begin{definition}
Let $\mathcal{M}$ be an abelian monoid with neutral element $e.$
\begin{itemize}
\item[(i)] The elements $m,n \in \mathcal{M}$ are associated if there is a unit $u \in \mathcal{M}^*$ such that $m=u\cdot n.$
\item[(ii)] If all divisors of the non-unit element $m$ are associated to either $e$ or $m,$ then $m$ is called irreducible.
\item[(iii)] If $ab=ac$ implies $b=c$ for all elements $b,c \in \mathcal{M},$ then the element $a$ is called cancellable.
\end{itemize}
 \end{definition}

\begin{proposition}[\cite{Bokor2021}]\label{prop1}
Every element in $\mathcal{M}_{n}^\mathrm{Diff}$ admits a connected sum decomposition into a homotopy sphere and irreducible manifolds. Moreover, 
all units of $\mathcal{M}_{n}^\mathrm{Diff}$ are homotopy spheres.
\end{proposition}

Not all the elements of the the monoid $\mathcal{M}_{n}^\mathrm{Diff}$ given in Proposition \ref{prop1} can be cancellable. For example, the manifold 
$\mathbb{S}^2 \times \mathbb{S}^{n-2}$ is not cancellable for every $n\geq 4,$ and thus $\mathcal{M}_{n}^\mathrm{Diff}$ is not a unique factorisation monoid. Then the question of which types of high-dimensional manifolds form a unique factorisation monoid has become important. In the following theorem, Smale and Wall answered this question and showed that in some special cases the monoid $\mathcal{M}_{n}^{\mathrm{Diff},\mathrm{hc}}$ is a unique factorisation monoid, \cite[Corollary 1.3]{smale} and \cite{Wall}.

\begin{theorem}[\cite{smale,Wall}]\label{SMALEWALL}
 The monoid $\mathcal{M}_{2n}^{\mathrm{Diff},\mathrm{hc}}$ is a unique factorisation monoid for $n\equiv 3,5,7 \; \mathrm{mod} \;8$ with $n\neq 15,31.$ Precisely, for any $W^{2n} \in \mathcal{M}_{2n}^{\mathrm{Diff},\mathrm{hc}}$, there is a decomposition
 $$W^{2n} = M^{2n}_1 \# M^{2n}_2 \# \ldots \# M^{2n}_{k}$$ 
 which is unique up to renumbering and rescaling the irreducible manifolds $M^{2n}_i$ by units.
  \end{theorem}

Reidemeister–Franz torsion is a classical topological invariant that encodes geometric and combinatorial information not detected by homology alone \cite{T. A. Chap,Franz,RCLC,Reidemeister-Franz}. Since its introduction, it has played an important role in manifold theory and related areas, most notably in the disproof of the Hauptvermutung \cite{Milnor2,Milnor}. Subsequent developments, particularly those due to Milnor, established a precise
relationship between Reidemeister–Franz torsion and algebraic invariants such as the Alexander polynomial, and clarified its behavior under fundamental topological constructions.

In recent years, there has been growing interest in the applied and computational aspects of topological invariants. Because of this, understanding how Reidemeister–Franz torsion behaves under decomposition and reconstruction of manifolds is of central importance. When a closed, oriented $2n$-dimensional manifold $W^{2n} \in \mathcal{M}_{2n}^{\mathrm{Diff},\mathrm{hc}}$ is obtained by gluing simpler manifolds $M_L^{2n}$ and $M_R^{2n}$, an explicit
multiplicative formula for the torsion allows one to reduce global computations to local data. Moreover, the torsion of $W^{2n}$ is the product of the torsions of $M_L^{2n},$ $M_R^{2n},$ and the torsion of $(2n-1)$-sphere $\mathbb{S}^{2n-1}$ times a corrective term $\mathbb{T}_{RF}\left(\mathcal{H}_{\ast}\right)$ coming from homologies of the chain complexes of the cell-decompositions of manifolds. Due to Milnor \cite[Theorem 3.1]{Milnor}, under the acyclicity assumption on the chain complexes of the cell-decompositions of such manifolds, the Reidemeister–Franz torsion is multiplicative up to the corrective term, which becomes trivial in the acyclic case.

There exist numerous explicit computations in low-dimensional settings \cite{EDEYS1,EDEYS2,Porti,SozOJM,SozMathScandinavia}. However, little is known about Reidemeister–Franz torsion for higher-dimensional manifolds, especially when the acyclicity condition is not satisfied. The lack of such computations limits the applicability of torsion, where acyclicity cannot always be
assumed and corrective terms must be carefully controlled.

The aim of this paper is to address this gap by presenting a class of higher-dimensional closed manifolds for which the multiplicativity of Reidemeister–Franz torsion holds without the assumption of acyclicity. We show that, for such manifolds, the corrective term in the multiplicative formula is trivial.

\begin{theorem}\label{theo1}
For $n\equiv 3,5,7 \; \mathrm{mod} \;8$ with $n\neq 15,31,$ let $W_p^{2n} \in \mathcal{M}_{2n}^{\mathrm{Diff},\mathrm{hc}}$ such that 
$$W_p^{2n} \cong M^{2n}_1 \# M^{2n}_2 \# \ldots \# M^{2n}_{p+1},$$ 
where the summands $M^{2n}_j \in \mathcal{M}_{2n}^{\mathrm{Diff},\mathrm{hc}}$ are irreducible $2n$-manifolds. Let $\mathbf{h}^{W_p^{2n}}_\nu,$ $\mathbf{h}^{\mathbb{S}_i^{2n-1}}_\eta,$ and $\mathbf{h}_0^{\overline{\mathbb{D}_i^{2n}}}=f^i_{\ast}(\varphi_0(\mathbf{c}_0))$ be respectively bases of $H_\nu(W_p^{2n}),$ $H_\eta(\mathbb{S}_i^{2n-1}),$ and $H_0(\overline{\mathbb{D}_i^{2n}})$ for 
$\nu\in \{0,\ldots,2n\},$ $\eta \in \{0,\ldots,2n-1\},$ $i \in \{1,\ldots,p\}.$ Then there is a basis $\mathbf{h}^{M^{2n}_{j}}_\nu$ of $H_\nu(M^{2n}_j)$ for each $j$ such that the following formula holds 
\begin{equation*}
|\mathbb{T}_{RF}(W_p^{2n},\{\mathbf{h}_\nu^{W_p^{2n}}\}_{0}^{2n})|=\prod_{j=1}^{p+1}|{\mathbb{T}_{RF}(M^{2n}_{j},\{\mathbf{h}_\nu^{M^{2n}_{j}}\}_{0}^{2n})}|.
\end{equation*}
Here, $f^i_{\ast}$ is the map induced by the simple homotopy equivalence 
$f^i:\{*\}\rightarrow \overline{\mathbb{D}_i^{2n}}$ and the map
$\varphi_0:Z_0(C_{\ast})\rightarrow H_0(C_{\ast})$ is the natural projection, and 
$\mathbf{c}_0$ denotes the geometric basis of $C_{0}(C_{\ast})$ in the chain complex $C_{\ast}(\{*\})$ of the point $\{*\}$.
\end{theorem}

\section{The Reidemeister-Franz torsion}
 \label{Sec:3}
 We give the required definitions and the basic facts about Reidemeister-Franz torsion and symplectic chain complex. Further information and the detailed proof can be found in \cite{Milnor,Porti,SozOJM,Witten} and the references therein.
\subsection{The Reidemeister-Franz torsion of a general chain complex}
\label{subsec1}

Assume that $V$ is a $k$-dimensional vector space over $\mathbb{R}$ and all bases of $V$ are ordered. Let $\mathbf{e}= (e_1,\ldots,e_k)$ and $\mathbf{f}=(f_1,\ldots,
f_k)$ be any bases of $V.$ Then the following equality holds
$$e_i=\sum_ {j=1}^k a_{ij} f_j, \;\;  i=1,\ldots,k,$$
where the transition matrix $(a_{ij})$ is invertible $(k \times k)$-matrix over 
$\mathbb{R}.$ We define the determinant of the transition matrix from basis 
$\mathbf{e}$ to basis $\mathbf{f}$ as
$$ [\mathbf{e}\rightarrow\mathbf{f}]=\det (a_{ij}) \in \mathbb{R}^{*}(=\mathbb{R}-\{0\})$$
with the following properties :
\begin{itemize}
\item[(i)] {$[\mathbf{e}\rightarrow\mathbf{e}] = 1,$}
\item[(ii)] {For a third basis $\mathbf{g}$ of $V,$ $[\mathbf{g}\rightarrow \mathbf{e}] = [\mathbf{g}\rightarrow \mathbf{f}] \cdot [\mathbf{f}\rightarrow\mathbf{e}],$}
\item[(iii)] {For $V=\{0\},$ $[0\rightarrow 0] = 1$ by using the convention $1\cdot 0=0.$}
\end{itemize}

Let $C_{\ast}$ be a chain complex of finite dimensional vector spaces over 
$\mathbb{R}$  
 $$C_{\ast}=(0\rightarrow C_n\stackrel{\partial_n}{\rightarrow} %
  C_{n-1}\rightarrow\cdots\rightarrow C_1 \stackrel{%
  \partial_1}{\rightarrow} C_0\rightarrow 0).$$
   For $p \in \{0,\ldots,n\},$  $H_p(C_{\ast})=Z_p(C_{\ast})/B_p(C_{\ast})$ denotes the $p$-th homology space of the chain complex $C_{\ast},$ where 
$$ B_p(C_{\ast})=\mathrm{Im}\{\partial_{p+1}:C_{p+1}\rightarrow C_{p}
\},$$
$$ Z_p(C_{\ast})=\mathrm{Ker}\{\partial_{p}:C_{p}\rightarrow
C_{p-1} \}.
$$
Consider the sequences with the inclusion $\imath$ and the natural projection $\varphi_p$.
\begin{equation}\label{Equation1}
0\rightarrow Z_p(C_\ast) \stackrel{\imath}{\hookrightarrow} C_p
\stackrel{\partial_p}{\rightarrow} B_{p-1}(C_\ast) \rightarrow 0,
\end{equation}
\begin{equation}\label{Equation2}
0\rightarrow B_p(C_\ast) \stackrel{\imath}{\hookrightarrow} Z_p(C_\ast)
\stackrel{\varphi_p}{\twoheadrightarrow} H_p(C_\ast) \rightarrow 0.
\end{equation}

The First Isomorphism Theorem says the sequence (\ref{Equation1}) is short exact and also the definition of $H_p(C_\ast)$ gives the short exactness of the sequence (\ref{Equation2}).Let $s_p:B_{p-1}(C_{\ast})\to C_p$ and 
 $\ell_p:H_p(C_{\ast})\to
Z_p(C_{\ast})$ be denote the section of  $\partial_p:C_p \to B_{p-1}(C_{\ast})$ and$\varphi_p:Z_p(C_{\ast})\to H_p(C_{\ast})$, respectively. Applying Spliting Lemma to the sequences (\ref{Equation1}) and (\ref{Equation2}), we can write the space $C_p$ as the direct sums of the spaces as follows
\begin{equation}\label{Equation0}
C_p=B_{p}(C_{\ast})\oplus \ell_p(H_p(C_{\ast }))\oplus s_p(B_{p-1}(C_{\ast})).
\end{equation}
For any bases $\mathbf{c_p}=\{c^1_{p},\ldots,c^{m_p}_p \}$, 
$\mathbf{b_p}=\{b^1_{p},\ldots,b^{m_p}_p \}$, 
and $\mathbf{h_p}=\{h^1_{p},\ldots,h^{n_p}_p\}$ of spaces
$C_p,$ $B_p(C_\ast),$ $H_p(C_\ast),$ if we consider equation (\ref{Equation0}), we can obtain a new basis for $C_p$ such as
$$\mathbf{b}_p\sqcup
\ell_p(\mathbf{h}_p)\sqcup s_p(\mathbf{b}_{p-1}).$$

Considering the above arguments, Milnor defined the Reidemeister-Franz torsion of a general chain complex as follows.
\begin{definition}(\cite{Milnor}) 
Reidemeister-Franz torsion of a general chain complex $C_{\ast}$ with respect to bases $\{\mathbf{c}_p\}_{0}^n,$ $\{\mathbf{h}_p\}_{0}^{n}$ is
defined by 
$$\mathbb{T}_{RF}\left( C_{\ast},\{\mathbf{c}_p\}_{0}^n,\{\mathbf{h}_p\}_{0}^n \right)
 =\prod_{p=0}^n \left[\mathbf{c}_p\rightarrow \mathbf{b}_p\sqcup \ell_p(\mathbf{h}_p)\sqcup
s_p(\mathbf{b}_{p-1})\right]^{(-1)^{(p+1)}}.$$ Here,
$[\mathbf{c}_p\rightarrow\mathbf{b}_p\sqcup \ell_p(\mathbf{h}_p)\sqcup
s_p(\mathbf{b}_{p-1})]$ is the determinant of
the transition matrix from the initial basis $\mathbf{c}_p$ to the obtained basis
 $\mathbf{b}_p\sqcup \ell_p(\mathbf{h}_p)\sqcup
s_p(\mathbf{b}_{p-1})$
of $C_{p}.$
\end{definition}

In \cite{Milnor}, Milnor showed that Reidemeister torsion is independent of the bases $\mathbf{b}_p,$ and sections $s_p,\ell_p$. But it depends on the bases 
$\mathbf{c}_p$ and $\mathbf{h}_p.$ Making a change $\mathbf{c}_p \mapsto \widetilde{\mathbf{c}}_p$ and $\mathbf{h}_p \mapsto \widetilde{\mathbf{h}}_p$ changes the Reidemeister-Franz torsion as follows
\begin{equation}\label{change-base-formula}
\mathbb{T}_{RF}\left( C_{\ast},\{\mathbf{c}'_p\}_{0}^n,\{\mathbf{h}'_p\}_{0}^n \right)=
\displaystyle\prod_{p=0}^n\left(\dfrac{[\mathbf{c}_p\rightarrow \mathbf{c}'_p]}
{[\mathbf{h}_p\rightarrow \mathbf{h}'_p]}\right)^{(-1)^p}
\mathbb{T}_{RF}\left( C_{\ast},\{\mathbf{c}_p\}_{0}^n,\{\mathbf{h}_p\}_{0}^n \right).
\end{equation}

Consider the short exact sequence of
chain complexes 
\begin{equation}\label{er123}
0\to A_{\ast }\stackrel{\imath}{\rightarrow} B_{\ast
}\stackrel{\pi}{\to} D_{\ast }\rightarrow 0
\end{equation}
and its long exact sequence 
\begin{eqnarray*}
&&\mathcal{H}_{\ast}: \cdots \longrightarrow  H_p(A_{\ast
})\overset{\imath_{p}^{\ast}}{\longrightarrow }H_p(B_{\ast})
\overset{_{\pi_{p}^{\ast}}}{\longrightarrow }H_p(D) \nonumber \\ 
&&\quad \quad \quad  \quad  \quad \quad\tikz\draw[->,rounded corners,nodes={asymmetrical rectangle}](5.30,0.4)--(5.30,0)--(1,0)--(1,-0.4)node[yshift=4.0ex,xshift=12.0ex] {$\delta_{p}$}; \nonumber\\ 
&&\quad \quad \quad H_{p-1}(A_{\ast
})\overset{\imath_{p-1}^{\ast}}{\longrightarrow }H_{p-1}(B_{\ast})
\overset{_{\pi_{p-1}^{\ast}}}{\longrightarrow }H_{p-1}(D) \nonumber\\ 
&&\quad \quad  \quad \quad \quad \quad\tikz\draw[->,rounded corners,nodes={asymmetrical rectangle}](5.60,0.4)--(5.60,0)--(1,0)--(1,-0.4)node[yshift=4.0ex,xshift=12.0ex] {$\delta_{p-1}$};\nonumber \\ 
&&\quad \quad\quad  H_{p-2}(A_{\ast
})\overset{\imath_{p-2}^{\ast}}{\longrightarrow }\cdots
 \end{eqnarray*}
Indeed, $\mathcal{H}_{\ast}$ is an exact (or acyclic) chain complex $C_{\ast}$ of length $3n+2$ with the spaces $C_{3p}(\mathcal{H}_{\ast})=H_p(D_{\ast}),$
$C_{3p+1}(\mathcal{H}_{\ast})=H_p(A_{\ast}),$ and
 $C_{3p+2}(\mathcal{H}_{\ast})=H_p(B_{\ast}).$ The bases $\mathbf{h}_p^D,$ $\mathbf{h}_p^A,$ and $\mathbf{h}_p^B$ are considered as bases for $C_{3p}(\mathcal{H}_{\ast}),$ $C_{3p+1}(\mathcal{H}_{\ast}),$ and $C_{3p+2}(\mathcal{H}_{\ast}),$
respectively. By using this set-up, Milnor gave the multiplicativity property of  Reidemeister-Franz torsion as follows.
\begin{theorem}[\cite{Milnor}]\label{MilA} Assume that
$\mathbf{c}^{A}_p,$ $\mathbf{c}^B_p,$  $\mathbf{c}^D_p,$
$\mathbf{h}^A_p,$ $\mathbf{h}^B_p,$ and $\mathbf{h}^D_p$ are respectively bases
of $A_p,$ $B_p,$ $D_p,$  $H_p(A_{\ast }),$ $H_p(B_{\ast }),$ and
$H_p(D_{\ast })$. Assume also that $\mathbf{c}^A_p,$
$\mathbf{c}^B_p,$ and $\mathbf{c}^D_p$ are compatible in the sense
that $[\mathbf{c}^A_p\sqcup
\widetilde{\mathbf{c}^D_p}\rightarrow \mathbf{c}^B_p]=\pm 1,$ where
$\pi_p\left(\widetilde{\mathbf{c}^D_p}\right)=\mathbf{c}^D_p.$ Then the following  formula holds
\begin{eqnarray*}
\mathbb{T}_{RF}(B_{\ast},\{\mathbf{c}^B_p\}_{0}^n,\{\mathbf{h}^B_p\}_{0}^n)&=&\mathbb{T}_{RF}(A_{\ast},\{\mathbf{c}^A_p\}_{0}^n,
\{\mathbf{h}^A_p\}_{0}^n) \; \mathbb{T}_{RF}(D_{\ast},\{\mathbf{c}^D_p\}_{0}^n,\{\mathbf{h}^D_p\}_{0}^n)\\
&& \times \;\mathbb{T}_{RF}(\mathcal{H}_{\ast},\{\mathbf{c}_{3p}\}_{0}^{3n+2},\{0\}_{0}^{3n+2}).
\end{eqnarray*}
\end{theorem}

\begin{definition}
The corrective term is the Reidemeister-Franz torsion of the long exact sequence $\mathcal{H}_{\ast}$ stated in Theorem~\ref{MilA} as 
$$\mathbb{T}_{RF}(\mathcal{H}_{\ast},\{\mathbf{c}_{3p}\}_{0}^{3n+2},\{0\}_{0}^{3n+2}).$$

\end{definition}
\begin{lemma}(\cite{EDErdal},\cite{Milnor})\label{crorectivetermlem}
Let $d$ be the dimension of the CW-complex $X.$ 
\begin{itemize}
\item [(i)]If all the chain complexes in (\ref{er123}) are acylic, then
\begin{equation*}
\mathbb{T}_{RF}\left( \mathcal{H}_{\ast},\{\mathbf{h}_p\}_{p=0}^{3d+2},\{0\}_{p=0}^{3d+2}\right)
 =1.
\end{equation*}
\item[(ii)] If at least one of the chain complex in (\ref{er123}) is non-acylic, then
\begin{equation*}
\mathbb{T}_{RF}\left( \mathcal{H}_{\ast},\{\mathbf{h}_p\}_{p=0}^{3d+2},\{0\}_{p=0}^{3d+2}\right)
 =\prod_{p=0}^{3d+2} \left[\mathbf{h}_p\rightarrow 
 \mathbf{h}'_p\right]^{(-1)^{(p+1)}}.
\end{equation*}
Here, $\mathbf{h}'_p=\mathbf{b}_p\sqcup s_p(\mathbf{b}_{p-1}).$ 
\end{itemize}
\end{lemma}

 Theorem~\ref{MilA} implies the following result. 
\begin{lemma}[\cite{SozMathScandinavia}]\label{sumlemma}
If $A_{\ast},$ $D_{\ast}$ are  chain complexes, and 
$\mathbf{c}^A_p,$ $\mathbf{c}^D_p,$ $\mathbf{h}^A_p,$ and
$\mathbf{h}^D_p$ are bases of $A_p,$ $D_p,$ $H_p(A_{\ast}),$ and
$H_p(D_{\ast}),$ respectively, then the following equality is valid
\begin{eqnarray}
\mathbb{T}_{RF}(A_{\ast}\oplus D_{\ast},\{ \mathbf{c}^A_p\sqcup
\mathbf{c}^D_p\}_0^n,\{ \mathbf{h}^A_p\sqcup
\mathbf{h}^D_p\}_0^n)&=&\mathbb{T}_{RF}(A_{\ast},\{ \mathbf{c}^A_p\}_0^n,\{
\mathbf{h}^A_p\}_0^n) \nonumber\\
&& \times \; \mathbb{T}_{RF}(D_{\ast},\{\mathbf{c}^D_p\}_0^n,\{
\mathbf{h}^D_p\}_0^n).
\end{eqnarray}
\end{lemma}

\subsection{Symplectic chain complex}
\label{subsec2}
 Witten introduced the notion of symplectic chain complex and then considering Reidemeister-Franz torsion for these complexes, he computed the volume of several moduli space of representations from the fundamental group of a Riemann surface to a compact gauge group \cite{Witten}. 
 
 Now we give the definition of the symplectic chain complex and the necessary results.
 
\begin{definition}A \emph{symplectic chain complex} $(C_{\ast},\partial_{\ast},\{\omega_{\ast,q-\ast}\})$ of length $q$ is a chain complex satisfies the following properties: 
\begin{itemize}
	\item[(1)] $q\equiv 2 \; ( \text{mod }4 ),$
	\item[(2)] For $p=0,\ldots,q/2$, there is a non-degenerate bilinear form 
	$$\omega_{p,q-p}:C_p\times C_{q-p}\to\mathbb{R}$$ such that
	\begin{itemize}
	\item[(i)] $\partial-$\textit{compatible}: $\omega_{p,q-p}(\partial_{p+1}a,b)=(-1)^{p+1}\omega_{p+1,q-(p+1)}(a,\partial_{n-p}b)$,
	\item[(i)] \textit{anti-symmetric}: $\omega_{p,q-p}(a,b)=(-1)^{p(q-p)}\omega_{q-p,p}(b,a).$
	\end{itemize}
\end{itemize}
\end{definition}
From $q\equiv 2 \;(\text{mod\ }4)$ it follows 
$\omega_{p,q-p}(a,b)=(-1)^{p}\omega_{q-p,p}(b,a).$ 
By using $\partial-$compatibility of the bilinear maps $\omega_{p,q-p}:C_p\times C_{q-p}\to \mathbb{R}$, they can be extend to homologies
\cite{SozOJM}.

\begin{definition} For a symplectic chain
complex $(C_{\ast
},\partial_{\ast},\{\omega_{\ast,q-\ast}\})$ of length $q$, the bases $\mathbf{c}_p$ and $\mathbf{c}_{q-p}$ of $C_p$ and $C_{q-p}$ are $\omega$-\emph{compatible} if the matrix of $\omega_{p,q-p}$ in bases $\mathbf{c}_p,$ $\mathbf{c}_{q-p}$ equals to 
\begin{displaymath}
 \left\{ \begin{array}{ll}
 {I}_{k\times k} & \textrm{, $p\ne q/2,$ }\\
 \left(
\begin{array}{cc}
  0_{l\times l} & {I}_{l\times l} \\
  {I}_{l\times l} & 0_{l\times l} \\
\end{array}
\right)  & \textrm{, $p=q/2.$}\\
\end{array} \right.
\end{displaymath}
Here, $k=\dim(C_p)=\dim(C_{q-p}),$ and 
$2l=\dim(C_{q/2}).$ 
\end{definition}
Every symplectic chain complex has $\omega-$compatible bases. So the existence of $\omega-$compatible bases enables to compute the Reidemeister-Franz torsion of an $\mathbb{R}-$symplectic chain complex. The reader is referred to \cite{Porti,Turaev} for more information.

\begin{theorem}[Theorem 3.0.15, {\cite{SozOJM}}]\label{Tor_Of_General_SymplecticA}
Let $(C_{\ast },\partial_{\ast},\{\omega_{\ast,q-\ast}\})$
be a symplectic chain complex with $\omega-$compatible bases. 
For each $p \in \{0,\ldots,q\}$ if $\mathbf{c}_p,$ $\mathbf{h}_p$ are any bases of $C_p,$ $H_p(C_{\ast}),$ respectively, then the formula is valid
$$\mathbb{T}_{RF}(C_{\ast},\{\mathbf{c}_p\}_{0}^q,\{\mathbf{h}_p\}_{0}^q)=
\prod_{p=0}^{(q/2)-1}\left(\det [\omega_{p,q-p}]\right)^{(-1)^p}
\sqrt{\det[\omega_{_{q/2,q/2}}]}^{(-1)^{q/2}}.$$ Here,
$\det[\omega_{p,q-p}]$ denotes the determinant of the matrix of the
non-degenerate pairing $[\omega_{p,q-p}]:H_p(C_{\ast})\times
H_{q-p}(C_{\ast})\to \mathbb{R}$ in the bases $\mathbf{h}_p,$
$\mathbf{h}_{q-p}.$
\end{theorem}

\subsection{The Reidemeister-Franz torsion of a manifold}
\label{subsec3}
Let $K$ be a cell decomposition of an $n$-dimensional manifold $M^n.$ Denote the set of $p$-cells by $C_p(K).$ Then $K$ canonically defines a chain complex $C_{\ast}(K)$ of free abelian groups as follows
 $$C_{\ast}(K)=(0\to
C_n(K)\stackrel{\partial_n}{\to }C_{n-1}(K)\to\cdots\to
        C_1(K)\stackrel{\partial_1}{\to} C_0(K)\to 0),$$ where $\partial_p$ is the boundary operator for $p\in \{1,\ldots,n\}.$ By orienting the $p$-cells and ordering $C_p(K),$ the chain complex $C_{\ast}(K)$ has a \emph{geometric basis} $\mathbf{c}_p=\{c_p^1,\cdots,c_{p}^{m_p}\}$ of $C_p(K)$ for each $p \in \{1,\ldots,n\}.$
 \begin{definition} (\cite{Milnor})
  Let $\mathbf{h}_p$ be a basis of $H_p(M^n)$ for $p \in \{1,\ldots,n\}$. Then the Reidemeister-Franz torsion of $M^n$ is defined by 
  $$\mathbb{T}_{RF}\left(C_{\ast}(K),\{\mathbf{c}_p\}_{0}^n,\{\mathbf{h}_p\}_{0}^n\right).$$ 
   \end{definition}   

 Following the arguments introduced in
\cite[Lemma~2.0.5]{SozOJM}, one can obtain the following lemma.
\begin{lemma}\label{lemcell}
 Reidemeister-Franz torsion of $M^n$ does not depend on the cell decomposition.
\end{lemma}

From Lemma \ref{lemcell}, we can conclude that the Reidemeister-Franz torsion of $M^n$ is well-defined. So we denote by $\mathbb{T}_{RF}(M^n,\{\mathbf{h}_p\}_{0}^n)$ the Reidemeister-Franz torsion of $M^n$ in the basis $\mathbf{h}_p$ of $H_p(M^n)$ for $p \in \{1,\ldots,n\}$.

\begin{theorem}[Theorem 0.1-Theorem 3.5, \cite{SozMathScandinavia}]\label{absltone}
 Let $M^n$ be an orientable closed connected $n$-dimensional manifold and let
$\mathbf{h}_p$ a basis of $H_p(M^n)$ for $p \in\{0,\ldots,n\}.$
\begin{itemize}
\item[(i)]{ if $n$ is odd, then \\
$$|\mathbb{T}_{RF}(M^n,\{\mathbf{h}_p\}_{0}^{n})|=1.$$}
\item[(ii)]{if $n$ is even, then
\begin{eqnarray*}
 && |\mathbb{T}_{RF}(M^n,\{\mathbf{h}_p\}_{0}^{n})|=\prod_{p=0}^{n/2-1}\left|\det \bigtriangleup^{M^n}_{p,n-p}(\mathbf{h}_p,\mathbf{h}_{n-p})\right|^{(-1)^p} \\
  & & \quad \quad\quad \quad\quad \quad\quad \quad\; \times \;\sqrt{|\det \bigtriangleup^{M^n}_{n/2,n/2}(\mathbf{h}_{n/2},\mathbf{h}_{n/2})|}^{(-1)^{n/2}}.
\end{eqnarray*} 
 }
\end{itemize}
Here, $\bigtriangleup^{M^n}_{p,n-p}(\mathbf{h}_p,\mathbf{h}_{n-p})$ indicates the matrix of intersection pairing 
$(\cdot,\cdot)_{p,n-p}:H_p(M^n)\times H_{n-p}(M^n)\rightarrow \mathbb{R}$ 
in bases $\mathbf{h}_p$, $\mathbf{h}_{n-p}.$
\end{theorem}

\begin{remark}\label{ornek1}
Let $\mathbf{h}_i^{\mathbb{S}^n}$ be the homology basis of the unit sphere 
$\mathbb{S}^n$ for each $i \in\{0,\ldots,n\}.$ By Theorem \ref{absltone}, we have
\begin{itemize}
\item[(i)]{ if $n$ is odd, then
$|T_{RF}(\mathbb{S}^n,\{\mathbf{h}_0^{\mathbb{S}^n},\mathbf{h}_n^{\mathbb{S}^n}\})|=1,$ }
\item[(ii)] { if $n$ is even, then $|T_{RF}(\mathbb{S}^n,\{\mathbf{h}_0^{\mathbb{S}^n},\mathbf{h}_n^{\mathbb{S}^n}\})|=|(\det \bigtriangleup^{\mathbb{S}^n}_{0,n}(\mathbf{h}_0^{\mathbb{S}^n}, \mathbf{h}_n^{\mathbb{S}^n}))|.$} 
\end{itemize}
\end{remark}

\section{Main results}
\label{sec3}
In the present paper, we consider the Reidemeister-Franz torsion with untwisted $\mathbb{R}$-coefficients. For a manifold $M^{n},$ we mean by $H_i(M^{n})$ the homology space $H_i(M^{n};\mathbb{R})$ with $\mathbb{R}$-coefficient. 
We denote by $\mathbb{D}^{2n}$  the open unit ball in  $\mathbb{R}^{2n}$ and by $\overline{\mathbb{D}^{2n}}$ the closed unit ball in $\mathbb{R}^{2n}.$ 
 
 As a warm up, we are going to start this section by calculating the torsion of the closed unit ball $\overline{\mathbb{D}^{2n}}.$ Next, to prove Theorem~\ref{theo1}, we need to calculate the Reidemeister-Franz torsion of  
 $W^{2n}=M_L^{2n}\# M_R^{2n}$ in terms of torsions of $M_L^{2n}-{\mathbb{D}^{2n}}$ and $M_R^{2n}-{\mathbb{D}^{2n}}$ (Theorem~\ref{theo2dimpimn4}). Later, we give a formula to calculate the torsion of ${M^{2n}}-{\mathbb{D}^{2n}}$ in terms of the torsion of ${M^{2n}}$ (Theorem~\ref{proes1}). Therefore, these calculations form a template for the homological algebraic calculations that we need for the proof of Theorem~\ref{theo1}. 

 \subsection{The Reidemeister-Franz torsion of a closed unit ball}
 \label{subsec4}
A closed unit ball $\overline{\mathbb{D}^{2n}}$ and any point $\{*\}$ are special complexes by \cite[Definition in Section 12.3]{Milnor}. Consider the simple homotopy equivalence $f:\{*\}\rightarrow \overline{\mathbb{D}^{2n}}$ together with \cite[Lemma 12.5]{Milnor}. Since the Reidemeister-Franz torsion is a simple homotopy invariant (due to Remark 2.8 (a) of the preprint by Porti: Reidemeister torsion, hyperbolic three-manifolds, and character varieties, 2016, arXiv:1511.00400), for the homology basis $\mathbf{h}^{\overline{\mathbb{D}^{2n}}}_0=f_{\ast}(\mathbf{h}^{\{*\}}_0)$ of $\overline{\mathbb{D}^{2n}}$ we obtain
\begin{equation}\label{diskandpt}
\mathbb{T}_{RF}(\overline{\mathbb{D}^{2n}},
\{\mathbf{h}^{\overline{\mathbb{D}^{2n}}}_0\})= \mathbb{T}_{RF}(\{*\},
\{\mathbf{h}^{\{*\}}_0\}).
\end{equation}

Let $K=\{e_0\}$ denote the single $0$-cell of $\{*\}.$  Consider the following chain complex
\begin{equation}\label{chaincomplex23}
  \begin{array}{ccc}
 C_{\ast}:=(0 \stackrel{\partial_1}{\rightarrow}  C_{0}(K) \stackrel{\partial_0}{\rightarrow} 0).
  \end{array}
\end{equation}
From the following equalities
\begin{eqnarray*}
&& B_0(C_{\ast})=\mathrm{Im}\{\partial_{1}:C_{1}( C_{\ast})\rightarrow C_{0}( C_{\ast})
\}=\{0\},\\
&& Z_0(C_{\ast})=\mathrm{Ker}\{\partial_{0}:C_{0}( C_{\ast})\rightarrow
C_{-1}( C_{\ast})\}=C_{0}(K),
\end{eqnarray*}
 it follows that the $0$-th homology of $\{*\}$ can be given as 
$$H_0(\{*\})=Z_0(C_{\ast})/B_0(C_{\ast})\cong C_{0}(K).$$ 
Then there are
the following short exact sequences
\begin{equation}\label{Equation1eks}
0\rightarrow Z_0(C_\ast) \stackrel{\imath}{\hookrightarrow} C_0(C_{\ast})
\stackrel{\partial_0}{\rightarrow} B_{-1}(C_\ast) \rightarrow 0,
\end{equation}
\begin{equation}\label{Equation2eks}
0\rightarrow B_0(C_\ast) \stackrel{\imath}{\hookrightarrow} Z_0(C_\ast)
\stackrel{\varphi_0}{\rightarrow} H_0(C_\ast) \rightarrow 0.
\end{equation}
Here, $\imath$ is the inclusion and $\varphi_0$ is the natural
projection. Assume that $s_0:B_{-1}(C_{\ast})\rightarrow C_0(C_{\ast})$ and
$\ell_0:H_0(C_{\ast})\rightarrow Z_0(C_{\ast})$ are sections of the homomorphisms
$\partial_0:C_0(C_{\ast}) \rightarrow B_{-1}(C_{\ast}),$
$\varphi_0:Z_0(C_{\ast})\rightarrow H_0(C_{\ast}),$
respectively. As $B_{0}(C_{\ast})=B_{-1}(C_{\ast})$ is trivial, the homomorphism 
$\varphi_0$ becomes an isomorphism. Hence, the section $\ell_0$ is the inverse of this isomorphism. By using sequences (\ref{Equation1eks}) and (\ref{Equation2eks}), we obtain 
\begin{equation}\label{klur2eks}
C_0(C_{\ast})= \ell_0 (H_0(C_{\ast})).
\end{equation}

Assume also that $\mathbf{h}^{\{*\}}_0$ is an arbitrary basis of 
$H_0(\{*\}).$ From equation (\ref{klur2eks}) it follows
\begin{equation}\label{d?sk1}
\mathbb{T}_{RF}(\{*\},\{\mathbf{h}_0^{\{*\}}\})=\left[\mathbf{c}_0 \rightarrow \ell_0(\mathbf{h}_0^{\{*\}}) \right].
\end{equation}

 Combining equations (\ref{diskandpt}) and (\ref{d?sk1}), we obtain the following result.
  \begin{proposition}\label{disktorsiyon1}
Let $\mathbf{h}_{0}^{{\overline{\mathbb{D}^{2n}}}}$ be a basis of 
$H_0({\overline{\mathbb{D}^{2n}}})$ which is the image of the basis 
$\mathbf{h}_{0}^{\{*\}}=\varphi_0(\mathbf{c}_0)$ of 
$H_0(\{*\})$ under $f_{\ast}.$ Then we have
\begin{equation*}
\mathbb{T}_{RF}({\overline{\mathbb{D}^{2n}}},
\{\mathbf{h}^{{\overline{\mathbb{D}^{2n}}}}_0\})=\left[\mathbf{c}_0\rightarrow \ell_0(\mathbf{h}_0^{\{*\}})\right]=\left[\mathbf{c}_0\rightarrow \ell_0(\varphi_0(\mathbf{c}_0))\right]=\left[\mathbf{c}_0\rightarrow \mathbf{c}_0\right]=1.
\end{equation*}
\end{proposition}

  \subsection{The Reidemeister-Franz torsion of $W^{2n}= M_L^{2n}\# M_R^{2n}$}
\label{subsec5}

\begin{lemma}\label{lemmahomolgy}
For any differentiable orientable closed $2n$-manifold $M^{2n}$, the homology space $H_{2n}({M^{2n}}-{\mathbb{D}}^{2n})$ is trivial.
\end{lemma}

\begin{proof}
Let us abuse the notation and denote the triangulations of respective manifolds by 
$\mathbb{S}^{2n-1}$, ${M^{2n}}-{\mathbb{D}}^{2n}$, and $W^{2n}= M^{2n}\# M^{2n}.$ There exists the natural short exact sequence of the chain complexes :
\begin{equation}\label{lemhom2}
0\to C_{\ast}(\mathbb{S}^{2n-1})\stackrel{\imath}{\rightarrow}
C_{\ast}({M^{2n}}-{\mathbb{D}}^{2n})\oplus C_{\ast}({M^{2n}}-{\mathbb{D}}^{2n})\stackrel{\pi}{\rightarrow} C_{\ast}(W^{2n})\to 0.
\end{equation}

Associated with the short exact sequence (\ref{lemhom2}), there is the Mayer-Vietoris long exact sequence

 \begin{eqnarray*}
\mathcal{H}_{\ast}: &&  0 \overset{\imath_{2n}^{\ast}} {\rightarrow}
 H_{2n}({M^{2n}}-{\mathbb{D}}^{2n})\oplus H_{2n}({M^{2n}}-{\mathbb{D}}^{2n})\overset{\pi_{2n}^{\ast}}{\rightarrow} \mathbb{R} \overset{\delta_{2n}}{\rightarrow} \mathbb{R} \\
 &&  \quad \quad   \quad \quad \quad  \quad  \quad     \quad   \quad \tikz\draw[->,rounded corners,looseness=0,nodes={asymmetrical rectangle}](5.79,0.4)--(5.79,0)--(1,0)--(1,-0.4) node[yshift=5.0ex,xshift=15.0ex] {$\imath_{2n-1}^{\ast}$}; \\ \nonumber
&& H_{2n-1}({M^{2n}}-{\mathbb{D}}^{2n})\oplus H_{2n-1}({M^{2n}}-{\mathbb{D}}^{2n}) \overset{\pi_{2n-1}^{\ast}}{\rightarrow} H_{2n-1}(W^{2n})\cdots \\ \nonumber
  \end{eqnarray*}  
  By the exactness of $\mathcal{H}_{\ast},$ we have 
  \begin{equation}\label{lemhom3}
  \mathbb{R} \cong \mathrm{Im}\,\delta_{2n}\oplus \mathrm{Ker}\,\delta_{2n} \cong \mathrm{Im}\,\delta_{2n}\oplus \mathrm{Im}\,\pi_{2n}^{\ast}.
  \end{equation}
 Assume $\mathrm{Im}\,\delta_{2n}=\{0\}$. By equation (\ref{lemhom3}), we have $\mathrm{Im}\,\pi_{2n}^{\ast}\cong \mathbb{R}.$ Since $\mathrm{Im}\,\imath_{2n}^{\ast}=\{0\}=\mathrm{Ker}\,\pi_{2n}^{\ast} $, we get the following contradiction on the dimensions of the vector spaces:
 
 \begin{eqnarray*}
   H_{2n}({M^{2n}}-{\mathbb{D}}^{2n})\oplus H_{2n}({M^{2n}}-{\mathbb{D}}^{2n})
   &\cong & \mathrm{Ker}\,\pi_{2n}^{\ast} \oplus \mathrm{Im}\,\pi_{2n}^{\ast} \\\nonumber
    &\cong &\{0\}  \oplus  \mathrm{Im}\,\pi_{2n}^{\ast} \\ \nonumber
     &\cong & \mathbb{R}.
 \end{eqnarray*}
Thus our assumption is wrong. Hence, $\mathrm{Im}\,\delta_{2n}\cong  \mathbb{R}$ and $\mathrm{Im}\,\pi_{2n}^{\ast}=\{0\}.$ From the fact that $\mathrm{Im}\,\imath_{2n}^{\ast}=\{0\}=\mathrm{Ker}\,\pi_{2n}^{\ast}$ it follows
 \begin{eqnarray*}
   H_{2n}({M^{2n}}-{\mathbb{D}}^{2n})\oplus H_{2n}({M^{2n}}-{\mathbb{D}}^{2n})
   &\cong & \mathrm{Ker}\,\pi_{2n}^{\ast} \oplus \mathrm{Im}\,\pi_{2n}^{\ast} \\\nonumber
    &\cong & \{0\}  \oplus \{0\}= \{0\}.
 \end{eqnarray*}
Therefore $ H_{2n}({M^{2n}}-{\mathbb{D}}^{2n})=\{0\}.$
\end{proof}

Let $W^{2n}$ be a $2$-fold connected sum of highly connected differentiable orientable closed $2n$-manifolds
 $$W^{2n}= M_L^{2n}\# M_R^{2n}.$$
Hence, we obtain the following desired result:
\begin{theorem}\label{theo2dimpimn4}
Let $\mathbf{h}^{W^{2n}}_\nu,$ $\mathbf{h}_{0}^{\mathbb{S}^{2n-1}},$ and $\mathbf{h}_{2n-1}^{\mathbb{S}^{2n-1}}=\delta_{2n}(\mathbf{h}_{2n}^{W^{2n}})$ be respectively bases of $H_\nu(W^{2n}),$ $H_{0}(\mathbb{S}^{2n-1}),$ and $H_{2n-1}(\mathbb{S}^{2n-1})$ for $\nu \in \{0,\ldots,2n\}.$ Then there exist bases $\mathbf{h}^{{M^{2n}_L}-{\mathbb{D}}^{2n}}_\nu$ and $\mathbf{h}^{{M_R^{2n}}-{\mathbb{D}}^{2n}}_\nu$ of $H_\nu({M_L^{2n}}-{\mathbb{D}}^{2n})$ and $H_\nu({M_R^{2n}}-{\mathbb{D}}^{2n})$ such that the corrective term 
becomes $1$ without the assumption of acyclicity and the following formula is valid
\begin{eqnarray*}
\mathbb{T}_{RF}(W^{2n},\{\mathbf{h}_\nu^{W^{2n}}\}_{0}^{2n})&=&\mathbb{T}_{RF}({M_L^{2n}}-{\mathbb{D}^{2n}},\{\mathbf{h}_\nu^{{M^{2n}_L}-{\mathbb{D}^{2n}}}\}_{0}^{2n})\\
&&\times \;\mathbb{T}_{RF}({M_R^{2n}}-{\mathbb{D}^{2n}},\{\mathbf{h}_\nu^{{M_R^{2n}}-{\mathbb{D}^{2n}}}\}_{0}^{2n}) \\
&& \times \;\mathbb{T}_{RF}(\mathbb{S}^{2n-1},\{\mathbf{h}_0^{\mathbb{S}^{2n-1}},0,\ldots,0,\mathbf{h}_{2n-1}^{\mathbb{S}^{2n-1}}\})^{-1}. 
\end{eqnarray*}
\end{theorem}

\begin{proof}
We abuse the notation and denote the triangulations of respective manifolds by 
$\mathbb{S}^{2n-1}$, ${M_L^{2n}}-{\mathbb{D}}^{2n}$, ${M_R^{2n}}-{\mathbb{D}}^{2n}$, and $W^{2n}.$ There exists the natural short exact sequence of the chain complexes :
\begin{equation}\label{shrtsq154}
0\to C_{\ast}(\mathbb{S}^{2n-1})\stackrel{\imath}{\rightarrow}
C_{\ast}({M_L^{2n}}-{\mathbb{D}}^{2n})\oplus C_{\ast}({M_R^{2n}}-{\mathbb{D}}^{2n})\stackrel{\pi}{\rightarrow} C_{\ast}(W^{2n})\to 0.
\end{equation}

Associated with the short exact sequence (\ref{shrtsq154}), there is the Mayer-Vietoris long exact sequence $\mathcal{H}_{\ast}$: 
\begin{small}
  \begin{eqnarray*}
 && 0\rightarrow  H_{2n}(\mathbb{S}^{2n-1}) \overset{\imath_{2n}^{\ast}} {\rightarrow} H_{2n}({M_L^{2n}}-{\mathbb{D}}^{2n})\oplus H_{2n}({M_R^{2n}}-{\mathbb{D}}^{2n})\overset{\pi_{2n}^{\ast}}{\rightarrow} H_{2n}(W^{2n})\\
 &&  \quad \quad \quad \quad \tikz\draw[->,rounded corners,looseness=0,nodes={asymmetrical rectangle}](9.30,0.4)--(9.30,0)--(1,0)--(1,-0.4) node[yshift=5.0ex,xshift=32.0ex] {$\delta_{2n}$}; \\ \nonumber
&& H_{2n-1}(\mathbb{S}^{2n-1})
 \overset{\imath_{2n-1}^{\ast}} {\rightarrow}H_{2n-1}({M_L^{2n}}-{\mathbb{D}}^{2n})\oplus H_{2n-1}({M_R^{2n}}-{\mathbb{D}}^{2n}) \overset{\pi_{2n-1}^{\ast}}{\rightarrow} H_{2n-1}(W^{2n}) \\ \nonumber
 && \quad \quad  \quad \quad \tikz\draw[->,rounded corners,looseness=0,nodes={asymmetrical rectangle}](10.60,0.4)--(10.60,0)--(1,0)--(1,-0.4) node[yshift=5.0ex,xshift=32.0ex] {$\delta_{2n-1}$}; \\ \nonumber
&& H_{2n-2}(\mathbb{S}^{2n-1})
 \overset{\imath_{2n-2}^{\ast}} {\rightarrow}H_{2n-2}({M_L^{2n}}-{\mathbb{D}}^{2n})\oplus H_{2n-2}({M_R^{2n}}-{\mathbb{D}}^{2n}) \overset{\pi_{2n-2}^{\ast}}{\rightarrow} H_{2n-2}(W^{2n})\nonumber\\
 &&\quad  \quad \quad \quad  \tikz\draw[->,rounded corners,looseness=0,nodes={asymmetrical rectangle}](10.60,0.4)--(10.60,0)--(1,0)--(1,-0.4) node[yshift=5.0ex,xshift=32.0ex] {$\delta_{2n-2}$}; \\ \nonumber
 && \quad \quad \quad \quad \quad  \quad \quad \quad \quad \; \; \;  \; \; \;  \; \; \; \; \; \; \cdots  \cdots\cdots \cdots \cdots  \cdots\cdots \cdots  \cdots\cdots \\ \nonumber
  &&  \quad \quad \quad  \tikz\draw[->,rounded corners,looseness=0,nodes={asymmetrical rectangle}](11.0,0.4)--(11.0,0)--(1,0)--(1,-0.4) node[yshift=5.0ex,xshift=34.0ex] {$\delta_{1}$}; \\ \nonumber
&& H_0(\mathbb{S}^{2n-1})
 \overset{\imath_{0}^{\ast}} {\rightarrow}H_{0}({M_L^{2n}}-{\mathbb{D}}^{2n})\oplus H_0({M_R^{2n}}-{\mathbb{D}}^{2n}) \overset{\pi_0^{\ast}}{\rightarrow} H_0(W^{2n})\overset{\delta_0}{\rightarrow} 0.\nonumber
  \end{eqnarray*}
  \end{small}  

 By Lemma \ref{crorectivetermlem}, the Reidemeister torsion of $\mathcal{H}_{\ast}$ satisfies the following formula
\begin{equation}
\mathbb{T}_{RF}(\mathcal{H}_{\ast},\{\mathbf{h}_p\}_{0}^{6n+2},\{0\}_{0}^{6n+2})
 =\prod_{0}^{6n+2} \left[\mathbf{h}_p\rightarrow 
 \mathbf{h}'_p\right]^{(-1)^{(p+1)}},
\end{equation}
where $\mathbf{h}'_p$ is the obtained new basis $\mathbf{b}_p\sqcup s_p(\mathbf{b}_{p-1})$ of $C_p(\mathcal{H}_{\ast})$ for all $p.$ As the Reidemeister-Franz torsion is independent of the bases $\mathbf{b}_p$ and sections $s_p,$ we can choose the appropriable bases $\mathbf{b}_p$ and sections $s_p$ to show that the existence of the bases $\mathbf{h}^{{M^{2n}_L}-{\mathbb{D}}^{2n}}_\nu$ and $\mathbf{h}^{{M_R^{2n}}-{\mathbb{D}}^{2n}}_\nu$ in which the corrective term $\mathbb{T}( \mathcal{H}_{\ast},\{\mathbf{h}_p\}_{0}^{6n+2},\{0\}_{0}^{6n+2})$ is equal to $1.$ 

Let $C_p(\mathcal{H}_{\ast})$ denote the vector spaces in $\mathcal{H}_{\ast}$ for $p \in \{0,\ldots,6n+2\}.$ By using the arguments given in Section \ref{Sec:3}, we have the following equation for each $p$
\begin{equation}\label{lngexcteq1}
C_p(\mathcal{H}_{\ast})=B_p(\mathcal{H}_{\ast})\oplus s_{_{p}}(B_{p-1}(\mathcal{H}_{\ast})).
\end{equation}

First we consider the following part of the long exact sequence $\mathcal{H}_{\ast}:$
  \begin{equation*}
0\overset{\delta_1}{\rightarrow} H_0(\mathbb{S}^{2n-1})
 \overset{\imath_{0}^{\ast}} {\rightarrow}H_{0}({M_L^{2n}}-{\mathbb{D}}^{2n})\oplus H_0({M_R^{2n}}-{\mathbb{D}}^{2n}) \overset{\pi_0^{\ast}}{\rightarrow} H_0(W^{2n})\overset{\delta_0}{\rightarrow} 0.\nonumber
  \end{equation*}
  By Hurewicz theorem, $H_1(W^{2n})\cong\pi_1(W^{2n})=0$. So, $\delta_1$ is a zero map. We use equation (\ref{lngexcteq1}) for the vector space $C_0(\mathcal{H}_{\ast})=H_0(W^{2n}).$ Since $\mathrm{Im}\,\delta_0$ is trivial, we get
\begin{equation}\label{ngexcteq1es}
C_0(\mathcal{H}_{\ast})=\mathrm{Im}\,\pi_0^{\ast}\oplus s_{_0}(\mathrm{Im}\,\delta_0)=\mathrm{Im}\,\pi_0^{\ast}.
\end{equation}
 Choosing the basis $\mathbf{h}^{\mathrm{Im}\,\pi_0^{\ast}}$ of $\mathrm{Im}\,\pi_0^{\ast}$ as $\mathbf{h}_0^{W^{2n}},$ we get that $\mathbf{h}_0^{W^{2n}}$ becomes the obtained basis $\mathbf{h}'_0$ of $C_0(\mathcal{H}_{\ast})$ by equation (\ref{ngexcteq1es}). As $\mathbf{h}_0^{W^{2n}}$ is also the initial basis 
 $\mathbf{h}_0$ of $C_0(\mathcal{H}_{\ast}),$ the following equation is valid
\begin{equation}\label{rsteq1}
[\mathbf{h}_0\rightarrow\mathbf{h}'_0]=1. 
\end{equation}

If we consider equation (\ref{lngexcteq1}) for the space $C_1(\mathcal{H}_{\ast})=H_{0}({M_L^{2n}}-{\mathbb{D}}^{2n})\oplus H_0({M_R^{2n}}-{\mathbb{D}}^{2n}),$ then we have
\begin{equation}\label{ngexcteq1es23} 
C_1(\mathcal{H}_{\ast})=\mathrm{Im}\,\imath_0^{\ast}\oplus s_{_{1}}(\mathrm{Im}\,\pi_0^{\ast}).
\end{equation}
In the previous step, the basis $\mathbf{h}^{\mathrm{Im}\,\pi_0^{\ast}}$ of 
$\mathrm{Im}\,\pi_0^{\ast}$ was chosen as $\mathbf{h}_0^{W^{2n}}.$ By the isomorphism between $\mathrm{Im}\,\imath_0^{\ast}$ and $H_0(\mathbb{S}^{2n-1}),$ we can take the basis $\mathbf{h}^{\mathrm{Im}\,\imath_0^{\ast}}$ of $\mathrm{Im}\,\imath_0^{\ast}$ as $\imath_0^{\ast}(\mathbf{h}_0^{\mathbb{S}^{2n-1}}).$ By using equation (\ref{ngexcteq1es23}), we can write the obtained basis of $C_1(\mathcal{H}_{\ast})$ as follows
$$\mathbf{h}'_1=\left\{\imath_0^{\ast}(\mathbf{h}_0^{\mathbb{S}^{2n-1}}),
s_{_{1}}(\mathbf{h}_0^{W^{2n}})\right\}.$$ 

As a reason of connectedness of the manifolds, $H_0({M^{2n}_L}-{\mathbb{D}^{2n}})$ and $H_0({M_R^{2n}}-{\mathbb{D}^{2n}})$ are one-dimensional subspaces of the $2$-dimensional space $C_1(\mathcal{H}_{\ast}).$ Thus, there are non-zero vectors $(a_{_{11}},a_{_{12}})$ and $(a_{_{21}},a_{_{22}})$ such that
\begin{eqnarray*}
& &\left\{a_{_{11}}\imath_0^{\ast}(\mathbf{h}_0^{\mathbb{S}^{2n-1}})+a_{_{12}}
s_{_{1}}(\mathbf{h}_0^{W^{2n}})\right\},\\
& &\left\{a_{_{21}}\imath_0^{\ast}(\mathbf{h}_0^{\mathbb{S}^{2n-1}})+a_{_{22}}
s_{_{1}}(\mathbf{h}_0^{W^{2n}})\right\}
\end{eqnarray*}
are bases of $H_0({M^{2n}_L}-{\mathbb{D}^{2n}})$ and $H_0({M_R^{2n}}-{\mathbb{D}^{2n}}),$ respectively. Indeed, the $2\times 2$ matrix $A=(a_{_{ij}})$ with entries in $\mathbb{R}$ is invertible. Let us take the bases of $H_0({M^{2n}_L}-{\mathbb{D}^{2n}})$ and $H_0({M_R^{2n}}-{\mathbb{D}^{2n}})$ as follows
\begin{eqnarray*}
&&\mathbf{h}_0^{{M^{2n}_L}-{\mathbb{D}^{2n}}}=\left\{(\det A)^{-1}[a_{_{11}}\imath_0^{\ast}(\mathbf{h}_0^{\mathbb{S}^{2n-1}})+a_{_{12}}
s_{_{1}}(\mathbf{h}_0^{W^{2n}})]\right\},\\
&& \mathbf{h}_0^{{M_R^{2n}}-{\mathbb{D}^{2n}}}=
 \left\{a_{_{21}}\imath_0^{\ast}(\mathbf{h}_0^{\mathbb{S}^{2n-1}})+a_{_{22}}
 s_{_{1}}(\mathbf{h}_0^{W^{2n}})\right\}.
\end{eqnarray*}
 Hence, $\mathbf{h}_1=\{\mathbf{h}_0^{{M^{2n}_L}-{\mathbb{D}^{2n}}}, \mathbf{h}_0^{{M_R^{2n}}-{\mathbb{D}^{2n}}}\}$ becomes the initial basis of $C_1(\mathcal{H}_{\ast})$ and we get
\begin{equation}\label{rsteq124}
[\mathbf{h}_1\rightarrow\mathbf{h}'_1]=1
\end{equation}

Considering the space $C_2(\mathcal{H}_{\ast})=H_0(\mathbb{S}^{2n-1})$ in equation (\ref{lngexcteq1}) and using the fact that 
$\mathrm{Im}\,\delta_1=\{0\},$ we can write the space $C_2(\mathcal{H}_{\ast})$ as follows
\begin{equation}\label{ngexcteq234} 
C_2(\mathcal{H}_{\ast})=\mathrm{Im}\,\delta_1\oplus s_{_{2}}(\mathrm{Im}\,\imath_0^{\ast})= s_{_{2}}(\mathrm{Im}\,\imath_0^{\ast}).
\end{equation}
 By equation (\ref{ngexcteq234}), $s_{_{2}}(\imath_0^{\ast}(\mathbf{h}_0^{\mathbb{S}^{2n-1}}))=\mathbf{h}_0^{\mathbb{S}^{2n-1}}$ becomes the obtained basis $\mathbf{h}'_2$ of $C_2(\mathcal{H}_{\ast}).$ Note that the initial basis 
 $\mathbf{h}_2$ of $C_2(\mathcal{H}_{\ast})$ is also $\mathbf{h}_0^{\mathbb{S}^{2n-1}}.$ So we obtain
\begin{equation}\label{rsteq17y}
[\mathbf{h}_2\rightarrow\mathbf{h}'_2]=1 
\end{equation}

For each $j \in \{1,\ldots,2n-2\}$, let us consider the following parts of $\mathcal{H}_{\ast}$:  
\begin{equation}\label{eqnnwy1}
H_j(\mathbb{S}^{2n-1}) \overset{\imath_j^{\ast}} {\rightarrow} H_{j}({M_L^{2n}}-{\mathbb{D}}^{2n})\oplus H_{j}({M_R^{2n}}-{\mathbb{D}}^{2n})
 \overset{\pi_i^{\ast}} {\rightarrow} H_{j}(W^{2n})\overset{\delta_j}{\rightarrow} H_{j-1}(\mathbb{S}^{2n-1}).
\end{equation}
 Now we denote the vector spaces (from right to left) in sequence (\ref{eqnnwy1}) as $C_{3j-1}(\mathcal{H}_{\ast})$, $C_{3j}(\mathcal{H}_{\ast}),$ $C_{3j+1}(\mathcal{H}_{\ast})$ and $C_{3j+2}(\mathcal{H}_{\ast}).$ Note that for $j=1,$ the spaces $C_{3}(\mathcal{H}_{\ast}), C_{4}(\mathcal{H}_{\ast}),C_{5}(\mathcal{H}_{\ast})$ are equal to $\{0\}$ and for $j \in \{2,\ldots,2n-2\}$ the spaces $C_{3j-1}(\mathcal{H}_{\ast})$ and $C_{3j+2}(\mathcal{H}_{\ast})$ are equal to $\{0\}.$ Using the convention $1\cdot0=1$ for each 
 $j \in \{2,\ldots,2n-2\},$ we get
\begin{eqnarray}\label{rsteqHK78}
&&[\mathbf{h}_{3}\rightarrow \mathbf{h}'_{3}] =1, \nonumber\\
&&[\mathbf{h}_{4}\rightarrow \mathbf{h}'_{4}] =1, \nonumber\\
&&[\mathbf{h}_{5}\rightarrow \mathbf{h}'_{5}] =1, \nonumber\\
&&[\mathbf{h}_{3j-1}\rightarrow \mathbf{h}'_{3j-1}] =1, \nonumber\\
&& [\mathbf{h}_{3j+2}\rightarrow \mathbf{h}'_{3j+2}] =1.
\end{eqnarray}

By the exactness of $\mathcal{H}_{\ast},$ we get the following isomorphism for each $j:$
\begin{equation}\label{iso1}
H_{j}({M_L^{2n}}-{\mathbb{D}}^{2n})\oplus H_{j}({M_R^{2n}}-{\mathbb{D}}^{2n})
\overset {\pi_i^{\ast}}{\cong} H_{j}(W^{2n}).
\end{equation}

If we use equation (\ref{lngexcteq1}) for $C_{{3j+1}}(\mathcal{H}_{\ast})=H_j(W^{2n}),$ then the triviality of $\mathrm{Im}\,\delta_j$ gives the following equality 
\begin{equation}\label{ngexcteq23e35} 
C_{{3j+1}}(\mathcal{H}_{\ast})=\mathrm{Im}\,\pi_j^{\ast}\oplus s_{_{3j+1}}(\mathrm{Im}\,\delta_j)=\mathrm{Im}\,\pi_j^{\ast}.
\end{equation}
Since ${\mathrm{Im}\,\pi_j^{\ast}}$ equals to $H_j(W^{2n}),$ we can take the basis $\mathbf{h}^{\mathrm{Im}\,\pi_j^{\ast}}$ of ${\mathrm{Im}\,\pi_j^{\ast}}$ as $\mathbf{h}_j^{W^{2n}}.$ By equation (\ref{ngexcteq23e35}), $\mathbf{h}_j^{W^{2n}}$ becomes the obtained basis $\mathbf{h}'_{{3j+1}}$ of $C_{{3j+1}}(\mathcal{H}_{\ast}).$ As the initial basis $\mathbf{h}_{{3j+1}}$ of $C_{{3j+1}}(\mathcal{H}_{\ast})$ is also $\mathbf{h}_j^{W^{2n}},$ we get
\begin{equation}\label{rsteESYT34}
 [\mathbf{h}_{{3j+1}}\rightarrow \mathbf{h}'_{{3j+1}}] =1. 
 \end{equation}
 
Let $C_{{3j+2}}(\mathcal{H}_{\ast})$ be $H_{j}({M_L^{2n}}-{\mathbb{D}}^{2n})\oplus H_{j}({M_R^{2n}}-{\mathbb{D}}^{2n}).$ Since $\mathrm{Im}\,\imath_j^{\ast}=\{0\},$ equation~(\ref{lngexcteq1}) turns into
\begin{equation}\label{ngexcteqew34}  
C_{{3j+2}}(\mathcal{H}_{\ast})=\mathrm{Im}\,\imath_j^{\ast}\oplus s_{_{3j+2}}(\mathrm{Im}\,\pi_j^{\ast})=s_{_{3j+2}}(\mathrm{Im}\,\pi_j^{\ast}).
\end{equation}
By the isomorphism 
$H_{j}({M_L^{2n}}-{\mathbb{D}}^{2n})\oplus H_{j}({M_R^{2n}}-{\mathbb{D}}^{2n}) \overset {\pi_j^{\ast}}{\cong} H_{j}(W^{2n}),$ the section $s_{_{3j+2}}$ can be considered as $(\pi_j^{\ast})^{-1}.$ In the previous step, the basis $\mathbf{h}^{\mathrm{Im}\,\pi_j^{\ast}}$ of ${\mathrm{Im}\,\pi_j^{\ast}}$  was chosen as 
$\mathbf{h}_j^{W^{2n}}.$ Equation (\ref{ngexcteqew34}) implies that $s_{_{3j+2}}(\mathbf{h}_j^{W^{2n}})$ is the obtained basis
$\mathbf{h}'_{3j+2}$ of $C_{{3j+2}}(\mathcal{H}_{\ast}).$ Recall that the given basis $\mathbf{h}_j^{W^{2n}}$ of $H_j(W^{2n})$ is 
$$\left\{\mathbf{h}_{j,1}^{W^{2n}},\ldots,\mathbf{h}_{j,\;{d_{1j}+d_{2j}}}^{W^{2n}}\right\},$$ where $(d_{1j}+d_{2j})$ is the rank of $H_j(W^{2n}).$ As
$H_j({M_L^{2n}}-{\mathbb{D}}^{2n})$ and $H_j({M_R^{2n}}-{\mathbb{D}}^{2n})$ are $d_{1j}$ and $d_{2j}$-dimensional subspaces of $(d_{1j}+d_{2j})$-dimensional space
$C_4(\mathcal{H}_{\ast}),$ respectively there are non-zero vectors 
$(b_{{\nu1}},\cdots,b_{{\nu\;d_{1 j}+d_{2 j}}}),$ $\nu \in \{1,\ldots,d_{1j}+d_{2j}\}$ such that
$$\left\{\sum_{i=1}^{d_{1j}+d_{2j}}b_{{\nu i}}s_{_{4}}(\mathbf{h}_{1,i}^{W^{2n}})\right\}_{\nu=1}^{d_{1j}} \; \mathrm{and} \; \left\{\sum_{i=1}^{d_{1j}+d_{2j}}b_{{\nu i}}s_{_{4}}(\mathbf{h}_{1,i}^{W^{2n}})\right\}_{\nu=d_{1j}+1}^{d_{1j}+d_{2j}}$$ are bases of $H_j({M^{2n}_L}-{\mathbb{D}^{2n}})$ and $H_j({M_R^{2n}}-{\mathbb{D}^{2n}}).$ Moreover, the $(d_{1j}+d_{2j})\times (d_{1j}+d_{2j})$ real matrix $B=(b_{_{\nu i}})$ is non-singular. Let the followings be respectively basis of the spaces $H_j({M^{2n}_L}-{\mathbb{D}^{2n}})$ and $H_j({M_R^{2n}}-{\mathbb{D}^{2n}})$ 
\begin{eqnarray*}
&&\mathbf{h}_j^{{M^{2n}_L}-{\mathbb{D}^{2n}}}=\left\{\det(B)^{-1}\sum_{i=1}^{d_{1j}+d_{2j}}b_{{1 i}}s_{_{4}}(\mathbf{h}_{1,i}^{W^{2n}}),\; \left\{\sum_{i=1}^{d_{1j}+d_{2j}}b_{{\nu i}}s_{_{4}}(\mathbf{h}_{1,i}^{W^{2n}})\right\}_{\nu=2}^{d_{1j}}\right\},\\
&& \mathbf{h}_j^{{M^{2n}_R}-{\mathbb{D}^{2n}}}=
\left\{\sum_{i=1}^{d_{1j}+d_{2j}}b_{{\nu i}}s_{_{4}}(\mathbf{h}_{1,i}^{W^{2n}})\right\}_{\nu=d_{1j}+1}^{d_{1j}+d_{2j}}.
\end{eqnarray*}
Choosing the initial basis $\mathbf{h}_{{3j+2}}$ of $C_{3j+2}(\mathcal{H}_{\ast})$ as 
$$\left\{\mathbf{h}_j^{{M^{2n}_L}-{\mathbb{D}^{2n}}},  \mathbf{h}_j^{{M^{2n}_R}-{\mathbb{D}^{2n}}}\right\},$$ we obtain
\begin{equation}\label{rsteq1ER7}
 [\mathbf{h}_{{3j+2}}\rightarrow \mathbf{h}'_{{3j+2}}] =1. 
 \end{equation}

 Now we consider the last part of the sequence $\mathcal{H}_{\ast}:$
\begin{small}
 \begin{eqnarray*}
 && 0\rightarrow H_{2n}(\mathbb{S}^{2n-1}) \overset{\imath_{2n}^{\ast}} {\rightarrow} H_{2n}({M_L^{2n}}-{\mathbb{D}}^{2n})\oplus H_{2n}({M_R^{2n}}-{\mathbb{D}}^{2n})\overset{\pi_{2n}^{\ast}}{\rightarrow} H_{2n}(W^{2n})\\
 &&\quad  \quad  \quad  \quad  \tikz\draw[->,rounded corners,looseness=0,nodes={asymmetrical rectangle}](9.30,0.4)--(9.30,0)--(1,0)--(1,-0.4) node[yshift=5.0ex,xshift=32.0ex] {$\delta_{2n}$}; \\
&& H_{2n-1}(\mathbb{S}^{2n-1})
 \overset{\imath_{2n-1}^{\ast}} {\rightarrow}H_{2n-1}({M_L^{2n}}-{\mathbb{D}}^{2n})\oplus H_{2n-1}({M_R^{2n}}-{\mathbb{D}}^{2n})\\
 &&\quad  \quad  \quad  \quad  \tikz\draw[->,rounded corners,looseness=0,nodes={asymmetrical rectangle}](5.45,0.4)--(5.45,0)--(1,0)--(1,-0.4) node[yshift=5.0ex,xshift=16.0ex] {${\pi_{2n-1}^{\ast}}$}; \\ 
 && H_{2n-1}(W^{2n}).
\end{eqnarray*}
\end{small}

If we use the convention $1\cdot0=1$ for the spaces $C_{6n+2}(\mathcal{H}_{\ast})=H_{2n}(\mathbb{S}^{2n-1})=\{0\}$ and $C_{6n+1}(\mathcal{H}_{\ast})=H_{2n}({M_L^{2n}}-{\mathbb{D}}^{2n})\oplus H_{2n}({M_R^{2n}}-{\mathbb{D}}^{2n})=\{0\}$, then we have
\begin{eqnarray}\label{rsteqlyt45}
&&[\mathbf{h}_{6n+1}\rightarrow \mathbf{h}'_{6n+1}] =1, \nonumber\\
&& [\mathbf{h}_{6n+2}\rightarrow \mathbf{h}'_{6n+2}] =1.
\end{eqnarray}
From the exactness of the sequence $\mathcal{H}_{\ast}$ it follows that $\imath_{2n-1}^{\ast}$ is a zero map. Hence, we have the following isomorphism 
$$H_{2n-1}({M_L^{2n}}-{\mathbb{D}}^{2n})\oplus H_{2n-1}({M_R^{2n}}-{\mathbb{D}}^{2n})\overset {\pi_{2n-1}^{\ast}}{\cong} H_{2n-1}(W^{2n}).$$
Now we consider the steps for the isomorphism in the equation (\ref{iso1}) and we apply them to the above isomorphism. For the spaces
$C_{{6n-3}}(\mathcal{H}_{\ast})=H_{2n-1}(W^{2n})$ and $C_{{6n-2}}(\mathcal{H}_{\ast})=H_{2n-1}({M_L^{2n}}-{\mathbb{D}}^{2n})\oplus H_{2n-1}({M_R^{2n}}-{\mathbb{D}}^{2n}),$ the following equalities hold
\begin{eqnarray}\label{rstystbl43}
&& [\mathbf{h}_{{6n-3}}\rightarrow\mathbf{h}'_{{6n-3}}] =1. \nonumber\\
&&  [\mathbf{h}_{{6n-2}}\rightarrow \mathbf{h}_{{6n-2}}] =1. 
\end{eqnarray}

Let us consider the space $C_{6n-1}(\mathcal{H}_{\ast})=H_{2n-1}(\mathbb{S}^{2n-1})$ in equation (\ref{lngexcteq1}). The equality 
$\mathrm{Im}\,  \imath_{2n-1}^{\ast}=\{0\}$ implies 
\begin{equation}\label{ngexcteqewty60b4} 
C_{6n-1}(\mathcal{H}_{\ast})=\mathrm{Im}\,\delta_{2n} \oplus s_{_{6n-1}}(\mathrm{Im}\,\imath_{2n-1}^{\ast})=\mathrm{Im}\,\delta_{2n}.
\end{equation}
 Recall that $\mathbf{h}_{2n-1}^{\mathbb{S}^{2n-1}}=\delta_{2n}(\mathbf{h}_{2n}^{W^{2n}})$ is the initial basis $\mathbf{h}_{6n-1}$ of $C_{6n-1}(\mathcal{H}_{\ast}).$ Taking the basis $\mathbf{h}^{\mathrm{Im}\, \delta_{2n}}$ of $\mathrm{Im}\, \delta_{2n}$ 
as $\mathbf{h}_{2n-1}^{\mathbb{S}^{2n-1}}$ 
and considering equation (\ref{ngexcteqewty60b4}) give that $\mathbf{h}_{2n-1}^{\mathbb{S}^{2n-1}}$ becomes the obtained basis $\mathbf{h}'_{6n-1}$ of $C_{6n-1}(\mathcal{H}_{\ast}).$ Hence, we get
\begin{equation}\label{rstyst67123rt}
[\mathbf{h}_{{6n-1}}\rightarrow \mathbf{h}'_{{6n-1}}] =1.
\end{equation}

Finally, let us consider equation (\ref{lngexcteq1}) for $C_{6n}(\mathcal{H}_{\ast})=H_{2n}(W^{2n}).$ By the fact that $\mathrm{Im}\,\pi_{2n}^{\ast}=\{0\},$ the following equality holds
\begin{equation}\label{ngexcteqewtstb4}
C_{{6n}}(\mathcal{H}_{\ast})=\mathrm{Im}\, \pi_{2n}^{\ast}\oplus s_{_{6n}}(\mathrm{Im}\,\delta_{2n})=s_{_{6n}}(\mathrm{Im}\,\delta_{2n}).
\end{equation}
 In the previous step, $\mathbf{h}_{2n-1}^{\mathbb{S}^{2n-1}}=\delta_{2n}(\mathbf{h}_{2n}^{W^{2n}})$ was chosen as the basis $\mathbf{h}^{\mathrm{Im}\,\delta_{2n}}$ of $\mathrm{Im}\,\delta_{2n}.$ From equation (\ref{ngexcteqewtstb4}) it follows that $s_{_{6n}}(\delta_{2n}(\mathbf{h}_{2n}^{W^{2n}}))=\mathbf{h}_{2n}^{W^{2n}}$ becomes the obtained basis 
$\mathbf{h}'_{6n}$ of $C_{6n}(\mathcal{H}_{\ast}).$ The initial basis $\mathbf{h}_{{6n}}$ of $C_{{6n}}(\mathcal{H}_{\ast})$ is also
$\mathbf{h}_{2n}^{W^{2n}},$ hence we get
\begin{equation}\label{etyrstysrt}
[\mathbf{h}_{{6n}}\rightarrow \mathbf{h}'_{{6n}}]=1. 
\end{equation}
If we consider equations 
(\ref{rsteq1}), (\ref{rsteq124}), (\ref{rsteq17y}), (\ref{rsteqHK78}), (\ref{rsteESYT34}), (\ref{rsteq1ER7}), (\ref{rsteqlyt45}), (\ref{rstystbl43}), (\ref{rstyst67123rt}), and (\ref{etyrstysrt}), then we show that the corrective term satisfies the following equation
\begin{equation}\label{lngrt45}
\mathbb{T}_{RF}( \mathcal{H}_{\ast},\{\mathbf{h}_p\}_{0}^{6n+2} ,\{0\}_{0}^{6n+2})
  =\prod_{p=0}^{6n+2} \left[\mathbf{h}_p\rightarrow 
  \mathbf{h}'_p\right]^{(-1)^{(p+1)}} =1. 
\end{equation}

The natural bases $\mathbf{c}_p^{W^{2n}},$ 
$\mathbf{c}_p^{{M_L^{2n}}-{\mathbb{D}^{2n}}},$  
$\mathbf{c}_p^{{M_R^{2n}}-{\mathbb{D}^{2n}}},$ and
$\mathbf{c}_p^{\mathbb{S}^{2n-1}}$ in the short exact sequence (\ref{shrtsq154}) are compatible. By combining Theorem~\ref{MilA}, Lemma~\ref{sumlemma}, and equation (\ref{lngrt45}), the following formula is valid
  \begin{eqnarray*}
\mathbb{T}_{RF}(W^{2n},\{\mathbf{h}_\nu^{W^{2n}}\}_{0}^{2n})&=&\mathbb{T}_{RF}({M_L^{2n}}-{\mathbb{D}^{2n}},\{\mathbf{h}_\nu^{{M^{2n}_L}-{\mathbb{D}^{2n}}}\}_{0}^{2n})\\
&& \times\;  \mathbb{T}_{RF}({M_R^{2n}}-{\mathbb{D}^{2n}},\{\mathbf{h}_\nu^{{M_R^{2n}}-{\mathbb{D}^{2n}}}\}_{0}^{2n})\\
&& \times \;\mathbb{T}_{RF}(\mathbb{S}^{2n-1},
\{\mathbf{h}_0^{\mathbb{S}^{2n-1}},0,\ldots,0,\mathbf{h}_{2n-1}^{\mathbb{S}^{2n-1}}\})^{-1}.
\end{eqnarray*}
\end{proof}

  \subsection{The Reidemeister-Franz torsion of ${M^{2n}}-{\mathbb{D}^{2n}}$}
  \label{subsec6}
  
\begin{theorem}\label{proes1}
Suppose that $M^{2n}$ is highly connected differentiable orientable closed $2n$-manifold. Then there is the following short exact sequence of the chain complexes 
\begin{equation*}
0\to C_{\ast}(\mathbb{S}^{2n-1})\stackrel{\imath}{\rightarrow} C_{\ast}({M^{2n} }-{\mathbb{D}^{2n}})\oplus C_{\ast}(\overline{\mathbb{D}^{2n}})\stackrel{\pi}{\rightarrow} C_{\ast}(M^{2n})\to 0
\end{equation*}
and its corresponding Mayer-Vietoris sequence 
 \begin{small}
  \begin{eqnarray*}
&\mathcal{H}_{\ast}:& 0 \rightarrow 
H_{2n}({M^{2n}}-{\mathbb{D}}^{2n})\overset{\pi_{2n}^{\ast}} {\longrightarrow} H_{2n}(M^{2n}) \overset{\delta_{2n}} {\longrightarrow} H_{2n-1}(\mathbb{S}^{2n-1})\\ \nonumber
&&\quad\quad \quad \quad \quad \quad  \quad \;  \tikz\draw[->,rounded corners,looseness=0,nodes={asymmetrical rectangle}](6.0,0.4)--(6.0,0)--(1,0)--(1,-0.4) node[yshift=4.5 ex,xshift=9.5ex] {$\imath_{2n-1}^{\ast}$}; \\ 
& &\quad \quad H_{2n-1}({M^{2n}}-{\mathbb{D}}^{2n})
 \overset{\pi_{2n-1}^{\ast}} {\longrightarrow} H_{2n-1}(M^{2n})\overset{\delta_{2n-1}}{\rightarrow} 0\\ \nonumber
 &&\quad\quad\quad\quad\;\;\; \tikz\draw[->,rounded corners,looseness=0,dashed](6.50,0.4)--(6.50,0)--(1,0)--(1,-0.4) node[yshift=4.5 ex,xshift=12.5ex] {}; \\
 & &\quad H_{2}({M^{2n}}-{\mathbb{D}}^{2n})
 \overset{\pi_{2}^{\ast}} {\longrightarrow} H_{2}(M^{2n})\overset{\delta_{2}}{\rightarrow} 0\\ \nonumber 
 && \quad \quad \quad \quad \quad\; \;\;\;\tikz\draw[->,rounded corners,looseness=0,nodes={asymmetrical rectangle}](4.10,0.4)--(4.10,0)--(1,0)--(1,-0.4) node[yshift=4.5 ex,xshift=11.0ex] {$\imath_{1}^{\ast}$}; \\ 
& &\quad\quad \; H_{1}({M^{2n}}-{\mathbb{D}}^{2n})
 \overset{\pi_{1}^{\ast}} {\longrightarrow} H_{1}(M^{2n}) \\ \nonumber 
&&\quad\quad \quad \quad \tikz\draw[->,rounded corners,looseness=0](4.10,0.4)--(4.10,0)--(1,0)--(1,-0.4) node[yshift=5.0ex,xshift=10.0ex] {$\delta_{1}$}; \\ 
  & &\quad \quad H_0(\mathbb{S}^{2n-1})
 \overset{\imath_{0}^{\ast}} {\longrightarrow} H_0({M^{2n} }-{\mathbb{D}^{2n}})\oplus H_0(\overline{\mathbb{D}^{2n}}) \overset{\pi_{0}^{\ast}}{\longrightarrow}H_0(M^{2n})\overset{\delta_0}{\longrightarrow} 0.
  \end{eqnarray*}
 \end{small}
 Suppose also that $\mathbf{h}^{{M^{2n}}-{\mathbb{D}^{2n}}}_\nu$ and 
 $\mathbf{h}^{\mathbb{S}^{2n-1}}_\eta$ are respectively bases of the homology spaces $H_\nu({M^{2n}}-{\mathbb{D}^{2n}}),$ $H_\eta(\mathbb{S}^{2n-1})$ for $\nu \in \{0,\ldots,2n\},$ $\eta \in \{0,\ldots,2n-1\},$ and $\mathbf{h}_0^{\overline{\mathbb{D}^{2n}}}=f_{\ast}(\varphi_0(\mathbf{c}_0))$ is the basis of $H_0(\overline{\mathbb{D}^{2n}}).$ Then there exists a basis $\mathbf{h}^{M^{2n}}_\nu$ of $H_\nu(M^{2n})$ such that the formula holds
\begin{eqnarray*}
\mathbb{T}_{RF}({M^{2n}}-{\mathbb{D}^{2n}},\{\mathbf{h}_\nu^{{M^{2n}}-{\mathbb{D}^{2n}}}\}_{0}^{2n})&=&
\mathbb{T}_{RF}(M^{2n},\{\mathbf{h}_\nu^{M^{2n}}\}_{0}^{2n}) \nonumber\\
&& \times \;\mathbb{T}_{RF}(\mathbb{S}^{2n-1},\{\mathbf{h}_\eta^{\mathbb{S}^{2n-1}}\}_{0}^{2n-1}).
\end{eqnarray*}
\end{theorem}

\begin{proof}

 For $j \in \{0,\ldots,2n\},$ consider the long exact sequence $\mathcal{H}_{\ast}$ as an exact complex
 $C_{\ast}$ of length $6n+2$ with 
 $$C_{3j}(\mathcal{H}_{\ast})=H_j(M^{2n}),$$
 $$C_{3j+1}(\mathcal{H}_{\ast})=H_j({M^{2n}}-{\mathbb{D}^{2n}})\oplus H_j(\overline{\mathbb{D}^{2n}}),$$
  $$C_{3j+2}(\mathcal{H}_{\ast})=H_j(\mathbb{S}^{2n-1})$$
  For each $j$, we use the following equation that is given in Section 
  \ref{Sec:3}:
\begin{equation}\label{lngeqalty1}
C_j(\mathcal{H}_{\ast})=B_j(\mathcal{H}_{\ast})\oplus s_{_{j}}(B_{j-1}(\mathcal{H}_{\ast})).
\end{equation}

By Hurewicz theorem, $H_1(M^{2n})\cong\pi_1(M^{2n})=0$. So, $\delta_1$ is a zero map. We first consider the following part of the 
long exact sequence $\mathcal{H}_{\ast}:$
  \begin{eqnarray*}
0 \overset{\delta_1} {\rightarrow}  H_0(\mathbb{S}^{2n-1})
 \overset{\imath_{0}^{\ast}} {\rightarrow} H_0({M^{2n} }-{\mathbb{D}^{2n}})\oplus H_0(\overline{\mathbb{D}^{2n}}) \overset{\pi_{0}^{\ast}}{\rightarrow}H_0(M^{2n})\overset{\delta_0}{\rightarrow} 0.
  \end{eqnarray*}
  
First, we use equation (\ref{lngeqalty1}) for the vector space $C_0(\mathcal{H}_{\ast})=H_0(M^{2n}).$ Since $\mathrm{Im}\,\delta_0$ is trivial, we get
\begin{equation}\label{Eqnlng1}
C_0(\mathcal{H}_{\ast})=\mathrm{Im}\,\pi_0^{\ast}\oplus s_{_0}(\mathrm{Im}\,\delta_0)=\mathrm{Im}\,\pi_0^{\ast}.
\end{equation}

As $\mathrm{Im}\,\pi_0^{\ast}$ is a one-dimensional space, there is a non-zero vector $(a_{_{11}},a_{_{12}})$ such that 
$$\mathbf{h}^{\mathrm{Im}\,\pi_0^{\ast}}=\left\{a_{_{11}}\pi_0^{\ast}(\mathbf{h}^{{M^{2n}}-{\mathbb{D}^{2n}}}_\nu)+a_{_{12}}\pi_0^{\ast}(\mathbf{h}_0^{\overline{\mathbb{D}^{2n}}})\right\}$$
 is the  basis of $\mathrm{Im}\,\pi_0^{\ast}.$ 
 From equation (\ref{Eqnlng1}) it follows that
$\mathbf{h}^{\mathrm{Im}\,\pi_0^{\ast}}$ is the obtained basis $\mathbf{h}'_0$ of $C_0(\mathcal{H}_{\ast}).$ If we choose the initial basis $\mathbf{h}_0$ (namely, $\mathbf{h}_0^{M^{2n}}$) of $C_0(\mathcal{H}_{\ast})$ as 
$\mathbf{h}^{\mathrm{Im}\,\pi_0^{\ast}},$ then we get
\begin{equation}\label{Eqnlng2}
[\mathbf{h}_0\rightarrow\mathbf{h}'_0]=1 
\end{equation}

Considering equation (\ref{lngeqalty1}) for $C_1(\mathcal{H}_{\ast})=H_0({M^{2n} }-{\mathbb{D}^{2n}})\oplus H_0(\overline{\mathbb{D}^{2n}}),$ the space $C_1(\mathcal{H}_{\ast})$ can be expressed as follows
\begin{equation}\label{Eqnlng3} 
C_1(\mathcal{H}_{\ast})=\mathrm{Im}\,\imath_0^{\ast}\oplus s_{_{1}}(\mathrm{Im}\,\pi_0^{\ast}).
\end{equation}
Recall that in the previous step we chose the basis of $\mathrm{Im}\,\pi_0^{\ast}$ as $\mathbf{h}^{\mathrm{Im}\,\pi_0^{\ast}}.$ Since $s_1$ is a section of 
$\pi_0^{\ast},$ the following equality holds
$$s_{1}(\mathbf{h}^{\mathrm{Im}\,\pi_0^{\ast}})= \{
a_{_{11}}\mathbf{h}^{{M^{2n}}-{\mathbb{D}^{2n}}}_0+a_{_{12}}\mathbf{h}_0^{\overline{\mathbb{D}^{2n}}}\}.$$ 
As $\mathrm{Im}\,\imath_0^{\ast}$ is a one-dimensional subspace of $C_1(\mathcal{H}_{\ast}),$ there is a non-zero vector 
 $(a_{_{21}},a_{_{22}})$ such that 
 $$\left\{a_{_{21}}\mathbf{h}^{{M^{2n}}-{\mathbb{D}^{2n}}}_0+a_{_{22}}\mathbf{h}_0^{\overline{\mathbb{D}^{2n}}}\right\}$$ is a basis of $\mathrm{Im}\,\imath_0^{\ast}$ and clearly $A=(a_{ij})$ is $(2\times 2)$-real matrix with non-zero determinant. If we take the basis of $\mathrm{Im}\,\imath_0^{\ast}$ as follows $$\mathbf{h}^{\mathrm{Im}\,\imath_0^{\ast}}=\left\{-(\det A)^{-1}\left[a_{_{21}}\mathbf{h}^{{M^{2n}}-{\mathbb{D}^{2n}}}_0+a_{_{22}}\mathbf{h}_0^{\overline{\mathbb{D}^{2n}}}\right]\right\},$$ then by equation (\ref{Eqnlng3}), 
$$\mathbf{h}'_1=\left\{\mathbf{h}^{\mathrm{Im}\,\imath_0^{\ast}},s_{1}(\mathbf{h}^{\mathrm{Im}\,\pi_0^{\ast}})\right\}$$ becomes the obtained basis of $C_1(\mathcal{H}_{\ast}).$ Since the initial basis of $C_1(\mathcal{H}_{\ast})$ is 
$$\mathbf{h}_1=\left\{\mathbf{h}^{{M^{2n}}-{\mathbb{D}^{2n}}}_0,\mathbf{h}_0^{\overline{\mathbb{D}^{2n}}}\right\},$$ the determinant of the transition matrix becomes $1;$ that is,
\begin{equation}\label{Eqnlng4}
  \left[\mathbf{h}_1\rightarrow \mathbf{h}'_1\right]=1.
\end{equation}

Next, let us consider the space $C_2(\mathcal{H}_{\ast})=H_0(\mathbb{S}^{2n-1})$ in equation (\ref{lngeqalty1}). Using the fact that $\mathrm{Im}(\delta_1)$ is a trivial space, we get
\begin{equation}\label{Eqnlng5}
C_2(\mathcal{H}_{\ast})=\mathrm{Im}(\delta_1)\oplus s_{2}(\mathrm{Im}(\imath_0^{\ast}))=s_{2}(\mathrm{Im}(\imath_0^{\ast})).
\end{equation}
 Recall that the basis $\mathbf{h}^{\mathrm{Im}\,\imath_0^{\ast}}$ of $\mathrm{Im}\,\imath_0^{\ast}$ was chosen as 
$$ \mathbf{h}^{\mathrm{Im}\,\imath_0^{\ast}}=\left\{-(\det A)^{-1}\left[a_{_{21}}\mathbf{h}^{{M^{2n}}-{\mathbb{D}^{2n}}}_0+a_{_{22}}\mathbf{h}_0^{\overline{\mathbb{D}^{2n}}}\right]\right\}$$ in the previous step. 
It follows from equation (\ref{Eqnlng5}) that $s_{2}(\mathbf{h}^{\mathrm{Im}\,\imath_0^{\ast}})$ is the obtained basis $\mathbf{h}'_2$ of $C_2(\mathcal{H}_{\ast}).$ If we take the initial basis
$\mathbf{h}_2$ (namely, $\mathbf{h}_0^{\mathbb{S}^{2n-1}}$) of  $C_2(\mathcal{H}_{\ast})$ as $s_{2}(\mathbf{h}^{\mathrm{Im}\,\imath_0^{\ast}})$, then we obtain
\begin{equation}\label{Eqnlng6}
 \left[\mathbf{h}_2\rightarrow \mathbf{h}'_2\right]=1.
\end{equation}

By the exactness of the sequence $\mathcal{H}_{\ast}$, Lemma \ref{lemmahomolgy}, and the First Isomorphism Theorem, we obtain the followings for each $i\in \{1,\ldots,2n-1\}$
\begin{itemize}
\item[(i)]{$\imath_{i}^{\ast}$ is zero map,} 
\item[(ii)]{$\delta_{i}$ is zero map,}
\item[(iii)]{$H_i({M^{2n} }-{\mathbb{D}^{2n}})\overset {\pi^{\ast}_{i}}{\cong} H_i(M^{2n}),$}
\item[(iv)]$H_{2n}(M^{2n}) \overset{\delta_{2n}} {\cong} H_{2n-1}(\mathbb{S}^{2n-1}).$
\end{itemize}
For each $i\in \{1,\ldots,2n-1\},$ by using the isomorphism 
$$H_i({M^{2n}}-{\mathbb{D}^{2n}})\overset {\pi^{\ast}_{i}}{\cong} H_i(M^{2n})$$ and given basis $\mathbf{h}^{{M^{2n}}-{\mathbb{D}^{2n}}}_i$ of $H_i({M^{2n}}-{\mathbb{D}^{2n}}),$ we can consider the basis of $\mathbf{h}^{M^{2n}}_i$ of $H_i(M^{2n})$ as $\pi^{\ast}_{i}(\mathbf{h}^{{M^{2n}}-{\mathbb{D}^{2n}}}_i).$ As $\imath_{i}^{\ast}$ is zero map, 
$\mathrm{Im}\,\imath_{i}^{\ast}=\{0\}$. Since ${\pi^{\ast}_{i}}$ is an isomorphism, its inverse can be considered as the section $s_{_i}$. As in the proof of 
Theorem~\ref{theo2dimpimn4}, we obtain 

\begin{equation}\label{lnlght2}
  \prod_{p=4}^{6n-1} \left[\mathbf{h}_p\rightarrow 
  \mathbf{h}'_p\right]^{(-1)^{(p+1)}} =1. 
\end{equation}

Now we consider the isomorphism $H_{2n}(M^{2n}) \overset{\delta_{2n}} {\cong} H_{2n-1}(\mathbb{S}^{2n-1})$ and given basis 
$\mathbf{h}^{\mathbb{S}^{2n-1}}_{2n-1}$ of $H_{2n-1}(\mathbb{S}^{2n-1}).$ 
Using the same arguments stated above, we take the basis 
$\mathbf{h}^{M^{2n}}_{2n}$ of $H_{2n}(M^{2n})$ as 
$\delta_{2n}^{-1}(\mathbf{h}^{\mathbb{S}^{2n-1}}_{2n-1}).$ Then we get 
\begin{equation}\label{lnlght34h}
  \prod_{p=6n}^{6n+2} \left[\mathbf{h}_p\rightarrow 
  \mathbf{h}'_p\right]^{(-1)^{(p+1)}} =1. 
\end{equation}

If we combine equations (\ref{Eqnlng2}), (\ref{Eqnlng4}), (\ref{Eqnlng6}), (\ref{lnlght2}), and (\ref{lnlght34h}), then we obtain  

\begin{equation}\label{Eqnlng7}
\mathbb{T}_{RF}(\mathcal{H}_{\ast},\{\mathbf{h}_p\}_{0}^{6n+2} ,\{0\}_{0}^{6n+2})=  \prod_{p=0}^{6n+2} \left[\mathbf{h}_p\rightarrow \mathbf{h}'_p\right]^{(-1)^{(p+1)}} =1. 
\end{equation}

Note that the natural bases $\mathbf{c}_p^{M^{2n}},$ 
$\mathbf{c}_p^{{M^{2n}}-{\mathbb{D}^{2n}}},$  
$\mathbf{c}_p^{\mathbb{S}^{2n-1}},$ and
$\mathbf{c}_p^{{\overline{\mathbb{D}^{2n}}}}$ in the short exact sequence (\ref{shrtsq154}) are compatible. From Theorem~\ref{MilA}, Lemma~\ref{sumlemma}, and equation (\ref{Eqnlng7}) it follows
\begin{eqnarray}\label{Eqnlng8}
\nonumber \mathbb{T}_{RF}({M^{2n}}-{\mathbb{D}^{2n}},\{\mathbf{h}_\nu^{{M^{2n}}-{\mathbb{D}^{2n}}}\}_{0}^{n})&=&
\mathbb{T}_{RF}(M^{2n},\{\mathbf{h}_\nu^{M^{2n}}\}_{0}^{2n}) \nonumber\\
&& \times \;\mathbb{T}_{RF}(\mathbb{S}^{2n-1},\{\mathbf{h}_\eta^{\mathbb{S}^{2n-1}}\}_{0}^{2n-1}) \nonumber\\
&& \times \;\mathbb{T}_{RF}(\overline{\mathbb{D}^{2n}},
\{\mathbf{h}_0^{\overline{\mathbb{D}^{2n}}}\})^{-1}.
\end{eqnarray}

Since $\mathbf{h}_0^{\overline{\mathbb{D}^{2n}}}=f_{\ast}(\varphi_0(\mathbf{c}_0))$ is the given basis of $H_0(\overline{\mathbb{D}^{2n}}),$  by Proposition~\ref{disktorsiyon1} we have 
\begin{equation}\label{Eqnlng9}
\mathbb{T}_{RF}({\overline{\mathbb{D}^{2n}}},
\{\mathbf{h}^{{\overline{\mathbb{D}^{2n}}}}_0\})=1.
\end{equation}

If we consider equations (\ref{Eqnlng8}) and (\ref{Eqnlng8}) together, we obtain the following formula. Hence, this finishes the proof of Proposition~\ref{proes1}:
\begin{eqnarray*}
\mathbb{T}_{RF}({M^{2n}}-{\mathbb{D}^{2n}},\{\mathbf{h}_\nu^{{M^{2n}}-{\mathbb{D}^{2n}}}\}_{0}^{n})&=&
\mathbb{T}_{RF}(M^{2n},\{\mathbf{h}_\nu^{M^{2n}}\}_{0}^{2n}) \nonumber\\
&& \times \;\mathbb{T}_{RF}(\mathbb{S}^{2n-1},\{\mathbf{h}_\eta^{\mathbb{S}^{2n-1}}\}_{0}^{2n-1}).
\end{eqnarray*}
\end{proof}

\subsection{The proof of Theorem \ref{theo1}}
\label{subsec7}
Assume that $W_p^{2n} \in \mathcal{M}_{2n}^{\mathrm{Diff},\mathrm{hc}}$, 
where $n\equiv 3,5,7 \; \mathrm{mod} \;8$ with $n\neq 15,31.$ Since the monoid 
$\mathcal{M}_{2n}^{\mathrm{Diff},\mathrm{hc}}$ is a unique factorisation monoid, there is a decomposition
 $$W_p^{2n} \cong M^{2n}_1 \# M^{2n}_2 \# \ldots \# M^{2n}_{p+1},$$ 
where the summands $M^{2n}_j \in \mathcal{M}_{2n}^{\mathrm{Diff},\mathrm{hc}}$ are irreducible $2n$-manifolds. Assume also that $\mathbf{h}^{W_p^{2n}}_\nu,$ $\mathbf{h}^{\mathbb{S}_i^{2n-1}}_\eta,$ and $\mathbf{h}_0^{\overline{\mathbb{D}_i^{2n}}}=f^i_{\ast}(\varphi_0(\mathbf{c}_0))$ are respectively bases of $H_\nu(W_p^{2n}),$ $H_\eta(\mathbb{S}_i^{2n-1}),$ and $H_0(\overline{\mathbb{D}_i^{2n}}),$ $\nu\in \{0,\ldots,2n\},$ $\eta \in \{0,\ldots,2n-1\},$ $i \in \{1,\ldots,p\},$ where $f^i_{\ast}$ is the map induced by the simple homotopy equivalence $f^i:\{*\}\rightarrow \overline{\mathbb{D}_i^{2n}}$ and 
$\varphi_0:Z_0(\{*\})\rightarrow H_0(\{*\})$ is the natural projection, and $\mathbf{c}^j_0$ is the geometric basis of $C_{0}(\{*\}).$ 


Under the above assumptions, we prove that there exists a basis $\mathbf{h}^{M^{2n}_{j}}_\nu$ of $H_\nu(M^{2n}_j)$ for each $j$ such that the Reidemeister-Franz torsion of $W_p^{2n}$ can be written as the product of the Reidemeister-Franz torsions of $M^{2n}_{j}$'s.

 For each $i \in \{1,\ldots,p\},$ let $W_i^{2n} \in \mathcal{M}_{2n}^{\mathrm{Diff},\mathrm{hc}}$ be an $(i+1)$-fold connected sum of orientable closed $2n$-manifolds; namely 
 $W_i^{2n}=\overset{i+1}{\underset {j=1}{\#}} M^{2n}_j,$ where $n\equiv 3,5,7 \; \mathrm{mod} \;8$ with $n\neq 15,31$. We consider $M_L^{2n}=W^{2n}_{i-1}$ and $M_R^{2n}=M^{2n}_{i+1}$ such that 
 $$W^{2n}_{i}=W^{2n}_{i-1}\#M^{2n}_{i+1}.$$ 
 Then we have the following short exact sequence
\begin{equation}\label{seqa1}
0\to C_{\ast}(\mathbb{S}_i^{2n-1})\rightarrow C_{\ast}(W^{2n}_{i-1}-{\mathbb{D}^{2n}_i})\oplus C_{\ast}(M^{2n}_{i+1}-{\mathbb{D}^{2n}_i})\rightarrow C_{\ast}(W^{2n}_{i})\to 0.
\end{equation}
Assume that $\delta_{2n}:H_{2n}(M_{i+1}^{2n})\rightarrow H_{2n-1}(\mathbb{S}_i^{2n-1})$ is a map in the long exact sequence associated to sequence (\ref{seqa1}) and  
$\mathbf{h}_{2n-1}^{\mathbb{S}^{2n-1}}=\delta_{2n}(\mathbf{h}_{2n}^{W^{2n}})$ is a basis of $H_{2n-1}(\mathbb{S}^{2n-1})$. By Theorem \ref{theo2dimpimn4}, for a given basis $\mathbf{h}^{W^{2n}_i}_\nu$ of $H_\nu(W^{2n}_i),$ there exist bases $\mathbf{h}^{{W_{i-1}^{2n}}-{\mathbb{D}_i^{2n}}}_\nu$ and $\mathbf{h}^{{{M_{i+1}^{2n}}-{\mathbb{D}_{i}^{2n}}}}_\nu$ of $H_\nu({W_{i-1}^{2n}}-{\mathbb{D}_i^{2n}})$ and $H_\nu({{M_{i+1}^{2n}}-{\mathbb{D}_{i}^{2n}}})$ such that the formula is valid
\begin{eqnarray}\label{seqa2}
\mathbb{T}_{RF}(W_i^{2n},\{\mathbf{h}_\nu^{W_i^{2n}}\}_{0}^{2n})&=&
\mathbb{T}_{RF}({W_{i-1}^{2n}}-{\mathbb{D}_i^{2n}},\{\mathbf{h}_\nu^{{W_{i-1}^{2n}}-{\mathbb{D}_i^{2n}}}\}_{0}^{2n}) \\
&& \times\; \mathbb{T}_{RF}({M_{i+1}^{2n}}-{\mathbb{D}_i^{2n}},\{\mathbf{h}_\nu^{{M_{i+1}^{2n}}-{\mathbb{D}_{i}^{2n}}}\}_{0}^{2n}) \nonumber\\
&& \times \;\mathbb{T}_{RF}(\mathbb{S}^{2n-1},
\{\mathbf{h}_0^{\mathbb{S}^{2n-1}},0,\ldots,0,\mathbf{h}_{2n-1}^{\mathbb{S}^{2n-1}}\})^{-1}\nonumber. 
\end{eqnarray}

Let us consider the following short exact sequences of chain complexes
\begin{equation}\label{seqa3}
0\to C_{\ast}(\mathbb{S}_i^{2n-1})\rightarrow C_{\ast}(W^{2n}_{i-1}-{\mathbb{D}^{2n}_i})\oplus C_{\ast}(\overline{\mathbb{D}_i^{2n}})\rightarrow C_{\ast}(W^{2n}_{i-1})\to 0,
\end{equation}
\begin{equation}\label{seqa3A}
 0\to C_{\ast}(\mathbb{S}_i^{2n-1})\rightarrow C_{\ast}(M^{2n}_{i+1}-{\mathbb{D}^{2n}_i})\oplus C_{\ast}(\overline{\mathbb{D}_i^{2n}})\rightarrow C_{\ast}(M^{2n}_{i+1})\to 0,
\end{equation}
and their associated Mayer-Vietoris long exact sequences as in 
Proposition~\ref{proes1}. In Proposition~\ref{proes1}, $\mathbf{h}_\nu^{W_{i-1}^{2n}}$ and 
$\mathbf{h}_\nu^{{M_{i+1}^{2n}}-{\mathbb{D}_{i}^{2n}}}$ are any given homology bases. So we can take these bases as above which is satisfying the equation (\ref{seqa2}). Since it is arbitrarily given basis in Proposition \ref{proes1}, we can respectively choose the same bases $\mathbf{h}_{0}^{\mathbb{S}^{2n-1}}$ and $\mathbf{h}_{2n-1}^{\mathbb{S}^{2n-1}}=\delta_{2n}(\mathbf{h}_{2n}^{W^{2n}})$ of $H_{0}(\mathbb{S}^{2n-1})$ and $H_{2n-1}(\mathbb{S}^{2n-1})$ for both sequences (\ref{seqa3}) and (\ref{seqa3A}). Hence, for the basis 
$\mathbf{h}_0^{\overline{\mathbb{D}_i^{2n}}}=f^i_{\ast}(\varphi_0(\mathbf{c}_0))$ of $H_0(\overline{\mathbb{D}_i^{2n}}),$ there exist respectively bases $\mathbf{h}_\nu^{W_{i-1}^{2n}}$ and $\mathbf{h}_\nu^{M_{i+1}^{2n}}$ of $H_\nu({W_{i-1}^{2n}})$ and $H_\nu({M_{i+1}^{2n}})$ such that the following formulas hold
 \begin{eqnarray}\label{seqa4}
 \mathbb{T}_{RF}({M_{i+1}^{2n}}-{\mathbb{D}^{2n}},\{\mathbf{h}_\nu^{{M_{i+1}^{2n}}-{\mathbb{D}^{2n}}}\}_{0}^{2n})&=&\mathbb{T}_{RF}(\mathbb{S}^{2n-1},
\{\mathbf{h}_0^{\mathbb{S}^{2n-1}},0,\ldots,0,\mathbf{h}_{2n-1}^{\mathbb{S}^{2n-1}}\})
 \nonumber\\
&& \times \;\mathbb{T}_{RF}(M_{i+1}^{2n},\{\mathbf{h}_\nu^{M_{i+1}^{2n}}\}_{0}^{2n}). 
\end{eqnarray}
\begin{eqnarray}\label{seqa5}
 \mathbb{T}_{RF}({W_{i-1}^{2n}}-{\mathbb{D}^{2n}},\{\mathbf{h}_\nu^{{W_{i-1}^{2n}}-{\mathbb{D}^{2n}}}\}_{0}^{2n})&=&\mathbb{T}_{RF}(\mathbb{S}^{2n-1},
\{\mathbf{h}_0^{\mathbb{S}^{2n-1}},0,\cdots,0,\mathbf{h}_{2n-1}^{\mathbb{S}^{2n-1}}\}) \nonumber\\
&&\times \;\mathbb{T}_{RF}(M_{i+1}^{2n},\{\mathbf{h}_\nu^{M_{i+1}^{2n}}\}_{0}^{2n}). 
\end{eqnarray}
By combining equations (\ref{seqa2}), (\ref{seqa4}), and (\ref{seqa5}), we obtain the Reidemeister-Franz torsion of $W_i^{2n}$ with untwisted $\mathbb{R}$-coefficients in these homology bases as follows 
\begin{eqnarray*}
\mathbb{T}_{RF}(W^{2n}_i,\{\mathbf{h}_\nu^{W^{2n}_i}\}_{0}^3)&=&\mathbb{T}_{RF}(W^{2n}_{i-1},\{\mathbf{h}_\nu^{W^{2n}_{i-1}}\}_{0}^3) \;
\mathbb{T}_{RF}(M^{2n}_{i+1},\{\mathbf{h}_\nu^{M^{2n}_{i+1}}\}_{0}^{2n})\\
&& \times \;\mathbb{T}_{RF}(\mathbb{S}^{2n-1},
\{\mathbf{h}_0^{\mathbb{S}^{2n-1}},0,\ldots,0,\mathbf{h}_{2n-1}^{\mathbb{S}^{2n-1}}\}).
\end{eqnarray*}
Let us follow the above arguments inductively. Then we have 
\begin{eqnarray}\label{seqa6}
\mathbb{T}_{RF}(W_p^{2n},\{\mathbf{h}_\nu^{W_p^{2n}}\}_{0}^{2n})&=& \prod_{i=1}^{p}{\mathbb{T}_{RF}(\mathbb{S}^{2n-1},
\{\mathbf{h}_0^{\mathbb{S}^{2n-1}},0,\ldots,0,\mathbf{h}_{2n-1}^{\mathbb{S}^{2n-1}}\})} \nonumber\\
&& \times\; \prod_{j=1}^{p+1}{\mathbb{T}_{RF}(M^{2n}_{j},\{\mathbf{h}_\nu^{M^{2n}_{j}}\}_{0}^{2n})}.
\end{eqnarray}
If we take the absolute value of both sides of equation (\ref{seqa6}), then by Theorem~\ref{absltone}(ii) we get 
\begin{equation*}
|\mathbb{T}_{RF}(W_p^{2n},\{\mathbf{h}_\nu^{W_p^{2n}}\}_{0}^{2n})|=\prod_{j=1}^{p+1}|{\mathbb{T}_{RF}(M^{2n}_{j},\{\mathbf{h}_\nu^{M^{2n}_{j}}\}_{0}^{2n})}|.
\end{equation*}
 This finishes the proof of Theorem \ref{theo1}.

\section{Declarations}
\label{Sec:5}
 The author has no competing interests to declare that are relevant to the content of this article.
\subsection{Ethical Approval}
 Research presented in this paper did not involve any sensitive data, therefore ethical approval is not necessary.
\subsection{Funding} 
The author was partially supported by TÜBİTAK under the project number 124F247.
\subsection{Availability of data and materials}
 Research presented in this paper has no associate data, so this declaration is not applicable.

\end{document}